\documentclass[11pt]{amsart}

\usepackage{amsmath, amsthm, amssymb}

\usepackage{lmodern, enumerate}
\usepackage{palatino}                       % text font
\usepackage[bigdelims,vvarbb]{newpxmath}    % math font
\usepackage[scaled=0.95]{inconsolata}       % Sans font
\usepackage[T1]{fontenc}                    % font encoding
\usepackage[normalem]{ulem}

\usepackage{comment}

\usepackage{tikz-cd}

\usepackage[abs]{overpic}		  % overpic automatically loads graphicx

\usepackage{xcolor}
\definecolor{indigo}{rgb}{0.29, 0.0, 0.51}  % custom colors
\usepackage[colorlinks, urlcolor=indigo, linkcolor=indigo, citecolor=indigo]{hyperref}

\usepackage[hcentering, top = 1.5in, total={5.9in, 8.2in}]{geometry}  % 1.4 vertical margin; 1.3 horizontal margin

\theoremstyle{plain}
\newtheorem{theorem}{Theorem}

\newtheorem{proposition}[theorem]{Proposition}
\newtheorem{lemma}[theorem]{Lemma}
\newtheorem{question}[theorem]{Question}

\newtheorem{mainthm}{Theorem}

\theoremstyle{definition}
\newtheorem{definition}[theorem]{Definition}

\theoremstyle{remark}
\newtheorem{remark}[theorem]{Remark}
\newtheorem{example}[theorem]{Example}
\newtheorem{exercise}[theorem]{Exercise}

\numberwithin{theorem}{section}

\newcommand{\dfn}[1]{{\bf #1}}        % definition
\newcommand{\R}{\mathbb{R}}           % the real numbers
\newcommand{\Z}{\mathbb{Z}}           % the integers
\newcommand{\C}{\mathbb{C}}           % the complex numbers
\newcommand{\N}{\mathbb{N}}           % the natural numbers
\DeclareMathOperator{\bd}{\partial}   % boundary
\makeatletter
\newcommand*\bigcdot{\mathpalette\bigcdot@{0.6}}
\newcommand*\bigcdot@[2]{\mathbin{\vcenter{\hbox{\scalebox{#2}{$\m@th#1\bullet$}}}}}
\makeatother

\DeclareMathOperator{\Diff}{Diff}         % group of diffeomorphisms
\DeclareMathOperator{\exot}{\text{exot}}
\DeclareMathOperator{\Exot}{\text{Exot}}
\DeclareMathOperator{\im}{\text{im}}

\newcommand{\Mod}{\operatorname{Mod}}
\newcommand{\FMod}{\operatorname{FMod}}

\newcommand{\crit}{\operatorname{crit}}
\newcommand{\Crit}{\operatorname{Crit}}

\newcommand{\bC}{\mathbb{C}}
\newcommand{\bD}{\mathbb{D}}

\newcommand{\bT}{\mathbb{T}}

\newcommand{\bZ}{\mathbb{Z}}
\newcommand{\cB}{\mathcal{B}}

\newcommand{\cU}{\mathcal{U}}

\usepackage{thmtools}
\usepackage{thm-restate}

\begin{document}

% title
\title{On nearly Lefschetz fibrations and spinal open books}

\author{Hyunki Min}

\author{Agniva Roy}

\author{Luya Wang}

\address{Department of Mathematics \\ University of Georgia \\ Athens, GA}
\email{hm32961@uga.edu}

\address{Department of Mathematics \\ Boston College \\ Chestnut Hill, MA}
\email{agniva.roy@bc.edu}

\address{Institute for Advanced Study \\ Princeton, NJ}
\email{luyawang@ias.edu}

%\subjclass[2020]{57R17}

% abstract
\begin{abstract}
We survey the recent development of spinal open books and nearly Lefschetz fibrations. The focus is on explaining past and new results regarding the detection of symplectic fillability and classification of symplectic fillings of contact 3-manifolds in light of the new machinery. We give many examples and a list of open questions. This can be considered as an additional chapter to the lecture notes by Ozbacgi \cite{ozbagci_survey} and Wendl \cite{Wendl_low_dim_curves}.
\end{abstract}

\maketitle
%\tableofcontents

%%%%%%%%%%%%%%%%%%%%%%%%
\section{Introduction}\label{sec: intro}
%%%%%%%%%%%%%%%%%%%%%%%%
\subsection{Overview and definitions} In this survey article, we provide a user’s guide to the series of papers \cite{spinal_I, spinal_II, spinal_exotic}, focusing on how to apply these results to classify and study symplectic fillings of contact $3$-manifolds supported by planar uniform spinal open books.

We first recall various types of symplectic fillings.

\begin{definition}
    A symplectic $4$-manifold $(W,\omega)$ is said to be a \textbf{weak symplectic filling} of the contact $3$-manifold $(M, \xi)$, if $M$ is the oriented boundary of $W$ and $\omega|_\xi > 0$.
\end{definition}

\begin{definition}
    A symplectic $4$-manifold $(W,\omega)$ is said to be a \textbf{strong symplectic filling} of the contact $3$-manifold $(M, \xi)$ if $M$ is the oriented boundary of $W$, and there exists a Liouville vector field $X$ near $\partial W$, i.e., $\mathcal{L}_X\omega = \omega$, and $\ker(\iota_X\omega) = \xi$.
\end{definition}
We sometimes refer to strong symplectic fillings as symplectic fillings.

\begin{definition}
    A symplectic $4$-manifold $(W,\omega)$ is said to be an \textbf{exact symplectic filling (or Liouville filling)} of the contact $3$-manifold $(M, \xi)$ if $M$ is the oriented boundary of $W$, there exists a global Liouville vector field $X$ on $W$, i.e., $\mathcal{L}_X\omega = \omega$, and $\ker(\iota_X\omega) = \xi$.
\end{definition}

\begin{definition}
    An exact filling  $(W,\omega)$ of $(M, \xi)$ is said to be a \textbf{Weinstein filling} if the Liouville vector field $X$ is gradient-like for a Morse function $f$ on $W$, i.e., $X \cdot\nabla f > 0$ away from the singularities of $X$, and these singularities coincide with critical points of $f$. By \cite{Cieliebak_Eliashberg} this is equivalent to a \textbf{Stein filling}. In light of this equivalence, in this note we will generally refer to both Weinstein and Stein fillings as Stein fillings.
\end{definition}

It is known that weak fillability implies tightness \cite{eli90}. Furthermore, the preceding definitions imply the following hierarchy among tight contact structures:
\[
    \text{Stein} \subsetneq \text{Exact} \subsetneq \text{Strong} \subsetneq \text{Weak} \subsetneq\text{Tight}
\]
It is also known that this hierarchy is proper. That is, there exist contact structures that are:
\begin{itemize}
    \item exactly fillable but not Stein fillable \cite{bowden}
    \item strongly fillable but not exactly fillable \cite{ghiggini,Min_Liouville}
    \item weakly fillable but not strongly fillable \cite{giroux_torsion, Eliashberg:T3,Gay:fillings}
    \item tight but not weakly fillable \cite{EH:nofillings,Christian:fillings}
\end{itemize}

However, in the case of contact $3$-manifolds supported by planar open book decompositions, a remarkable simplification occurs: Wendl \cite{Wendl_strongly} proves that all minimal (weak, strong, and exact) fillings are symplectically deformation equivalent to Stein fillings. Furthermore, these fillings are supported by Lefschetz fibrations over a disk.

\begin{theorem}[Wendl \cite{Wendl_strongly}]\label{thm: wendl_strongly}
    Suppose a contact $3$-manifold $(M,\xi)$ is supported by a planar open book $(P, \phi)$. Then, any minimal symplectic filling of $(M,\xi)$ can be foliated by $J$-holomorphic curves that are diffeomorphic to $P$ and admits a positive allowable Lefschetz fibration structure that induces the open book $(P,\phi)$ on its boundary.

    Moreover, there is a one-to-one correspondence between the deformation equivalence classes of the fillings and the set of positive allowable Lefschetz fibrations up to fiber-preserving diffeomorphism.
\end{theorem}

\begin{definition}
    Let $P$ be a compact oriented surface and $\phi \in \Mod_{\partial}(P)$. A \textbf{positive admissible factorization} of $\phi$ is a representative written as a product of positive Dehn twists along homologically essential curves on $P$.
\end{definition}

Combined with the fact that two Hurwitz equivalent positive admissible factorizations of the monodromy induce equivalent symplectic Lefschetz fibrations, Theorem~\ref{thm: wendl_strongly} {\em almost} reduces the classification of symplectic filling to a monodromy factorization problem.

\begin{remark}
    Theorem~\ref{thm: wendl_strongly} does not imply that there is a one-to-one correspondence between the set of minimal symplectic fillings and the set of positive admissible factorizations. Two different factorizations may induce equivalent symplectic fillings.
\end{remark}

Theorem~\ref{thm: wendl_strongly} applies to a large class of contact $3$-manifolds, which are supported by planar open books -- simply referred to as \textbf{planar contact $3$-manifolds}. This monodromy factorization technique has been used to classify the symplectic fillings of some contact structures on lens spaces \cite{Plamenevskaya_VHM} and other large classes of planar contact $3$-manifolds \cite{Kaloti:fillings, kaloti-li}. However, since the monodromy factorization problem is difficult in general, the classification of symplectic fillings for lens spaces, in full generality, has been achieved using different techniques. For example, for universally tight contact structures on lens spaces, Lisca \cite{Lisca:lensFillings} used the complements of divisor configurations in rational surfaces, and in the virtually overtwisted case, Christian--Li and Etnyre--Roy \cite{christian_li, etnyre_roy} used the ``mixed torus'' technique \cite{christian-menke}.

A major limitation of Theorem~\ref{thm: wendl_strongly} is that there exist contact $3$-manifolds that cannot be supported by planar open books. The first obstruction was found by Etnyre \cite{Etnyre_planar_obd}.

\begin{theorem}[{\cite[Theorem 4.1]{etnyre_obd}}] \label{thm: etnyre_negativedefinite}
    Let $(Y,\xi)$ be a contact $3$-manifold supported by a planar open book, then any symplectic filling of $(Y, \xi)$ is negative definite.
\end{theorem}

According to this theorem, the unique Stein fillable contact structure on $\bT^3$ does not admit a planar open book decomposition, since its Stein filling $\bT^2 \times \bD^2$, the disk cotangent bundle of the torus, is indefinite.

Nevertheless, Wendl \cite{Wendl_strongly} classifies the symplectic fillings of this contact structure on $\bT^3$ using techniques similar to the proof of Theorem \ref{thm: wendl_strongly}, by utilizing a topological decomposition similar to a planar open book. This decomposition can now be interpreted as a planar spinal open book (see Definition~\ref{def:spinalobd}), and the aforementioned classification result can be recovered using the notion of positive allowable monodromy factorizations. We discuss this example in detail in
\textsection \ref{subsubsec: T_sobd} and \textsection \ref{subsec: T3}.

The primary goal of this article is to explain how symplectic fillings of contact $3$-manifolds supported by \textbf{planar uniform spinal open books} can be classified using monodromy factorizations. This technique is applicable to a much larger class of contact manifolds than those supported by planar open books. The main results are the following theorems from \cite{spinal_I,spinal_II,spinal_exotic}.

\begin{mainthm} \label{thm: main}
    Let $(M, \xi)$ be a contact $3$-manifold supported by a planar uniform spinal open book and $(W,\omega)$ a minimal strong symplectic filling of $(M,\xi)$. Then, $(W,\omega)$ is symplectic deformation to a \textbf{positive allowable nearly Lefschetz fibration}, i.e., the complement of a neighborhood of positive multisections in a bordered Lefschetz fibration.
\end{mainthm}

Theorem \ref{thm: main} can be interpreted in terms of monodromy factorizations as follows.

\begin{mainthm}\label{thm: monodromy_fact}
    Let $(M, \xi)$ be a contact 3-manifold supported by a planar spinal open book $\pi$ and $\mathcal{B}$ a collection of surfaces that $\pi$ is uniform with respect to. Then a minimal strong filling of $(M, \xi)$ corresponds to a positive admissible factorization of the monodromy of the spinal open book with respect to some surface $B \in \mathcal{B}$.
\end{mainthm}

In \S\ref{sec: examples}, we outline a strategy for the classification of symplectic fillings based on Theorem~\ref{thm: monodromy_fact} and provide several examples.

\begin{remark}
    Note that Theorem~\ref{thm: monodromy_fact} does not claim a one-to-one correspondence. There may be redundancies in the correspondence, and we will see such examples in \S\ref{sec: examples}.
\end{remark}

In upcoming work \cite{breen-roy-wang_panlf}, it is shown that a positive admissible nearly Lefschetz fibration supports a Stein structure, generalizing the result on positive allowable Lefschetz fibrations (PALFs) by \cite{Akbulut_Ozbagci, spinal_I}\footnote{The ``positive'' part of the definition can be omitted, since all Lefschetz singularities are by definition positive. We keep it here as a historical continuation of terminology}. In particular, all minimal (weak, strong, exact) symplectic fillings of planar uniform spinal open books are symplectically deformation equivalent to Stein fillings. Details about the definition of nearly Lefschetz fibrations, and an associated Stein structure, will be given in \textsection \ref{sec: PANLF}.

\begin{theorem} \label{thm: Stein}
    A positive admissible nearly Lefschetz fibration (PANLF) supports a canonical Stein structure.
\end{theorem}

We will define the term \textbf{uniform} condition in \textsection \ref{sec: spinal_obd}, but for now, one can understand that this is some natural (technical) condition imposed on a planar spinal open book bounding a nearly Lefschetz fibration. The precise notion of a \textbf{positive admissible factorization} is clarified in Definition \ref{def: pos_adm_fact}.

\subsection{Organization} In \textsection\ref{sec: spinal_obd} we give definitions pertaining to spinal open books, and also discuss examples and methods of obtaining spinal open book decompositions. In \textsection\ref{sec: PANLF} we define and discuss nearly Lefschetz fibrations, and describe associated Stein structures. Then in \textsection\ref{sec: examples} we show how Theorems~\ref{thm: main} and \ref{thm: monodromy_fact} can be used to classify symplectic fillings for some families of $3$-manifolds. We conclude with some questions for future research in \textsection\ref{sec: questions}.

\subsection{Acknowledgements}
AR acknowledges support from an AMS-Simons travel grant. LW acknowledges support from NSF Grant DMS-2303437. We are grateful to the anonymous referee for their suggestions and corrections.

\subsection{Statement on AI use} No generative AI was used in the preparation of this document.

%%%%%%%%%%%%%%%%%%%%%%%%%%%%%%%%%%%%%%%%%%%%%%%%
\section{Spinal open books}\label{sec: spinal_obd}
%%%%%%%%%%%%%%%%%%%%%%%%%%%%%%%%%%%%%%%%%%%%%%%%

In this section, we review basic definitions and tools used to study spinal open book decompositions.

%----------------------------------------------------------------------------
\subsection{Spinal open book decompositions}\label{subsec: spinalobd}
%----------------------------------------------------------------------------
In this subsection, we review the definition of spinal open book decompositions, which first appeared in \cite{spinal_I, spinal_II}. We then recall the contact structures they support, and characterize the spinal open books that bound (nearly) Lefschetz fibrations.

\begin{definition}\label{def:spinalobd}
    A \textbf{spinal open book decomposition} of a closed oriented $3$-manifold $M$ is a decomposition $M = M_\Sigma \cup M_P$ (called the \textbf{spine} and \textbf{paper}, respectively) together with a pair of fibrations $\pi = (\pi_{\Sigma}, \pi_{P})$ such that
    \begin{align*}
        \pi_\Sigma &: M_\Sigma \to \Sigma \\
        \pi_P &: M_P \to S^1
    \end{align*}
    where
    \begin{enumerate}
        \item $\Sigma$ is a compact oriented surface whose connected components (called \textbf{vertebrae}) have nonempty boundary and $\pi_\Sigma$ is a trivial fibration with $S^1$ fiber;
        \item The fibers of $\pi_P$ are compact oriented surfaces whose connected components (called \textbf{pages} and denoted by $P$) have nonempty boundary and intersect transversely to $\partial M_P$. The intersection of $P$ and $M_\Sigma$ consists of fibers of $\pi_\Sigma$;
        \item For every component of $\partial M_P$, there exist local coordinates $(\phi, t, \theta) \in S^1 \times (-1,0] \times S^1$ on a collar neighborhood such that $\pi_P(\phi,t,\theta) = m\phi$ for some $m \in \mathbb{N}$. On every component of $\partial M_\Sigma$, there exist local coordinates $(s, \phi, \theta) \in (-1,0] \times S^1 \times S^1$ such that $\pi_{\Sigma}(s, \phi, \theta) = (s,\phi)$. The number $m$ is called the {\bf multiplicity} of $\pi_P$ at that boundary component of $M_P$. On the overlap between $M_P$ and $M_\Sigma$ the 2-torus coordinates $(\phi, \theta)$ agree, while the interval coordinates are related by $s = -t$;
        \item Fix a trivialization of the spine, i.e., an identification with $\Sigma \times S^1$, such that given $x \in \partial \Sigma$, $\{x\}\times S^1$ is the boundary component of a page. Then the paper $M_P$ can be identified with a mapping torus
        \[
            M_P = \R \times P \big/ \sim, \quad ( \tau, p) \sim (\tau-1,\phi(p))
        \]
        and $\pi_P: M_P \to S^1$ is given by $[(\tau, p)] \to [\tau]$. Here, $\phi \in \FMod(P)$ is the monodromy in the framed mapping class group of $P$. See \S\ref{subsec: framed_mcg} for the definition.
    \end{enumerate}
\end{definition}

Note that when all vertebrae are disks and the multiplicity of each component of $\partial M_P$ is $1$, the definition of spinal open books recovers the definition of ordinary open books. A contact form $\lambda$ on a spinal open book is called a \textbf{Giroux form} if $d\lambda$ is positive on the interior of each page and the Reeb vector field associated to $\lambda$ is positively tangent to every fiber of $\pi_\Sigma$. We say that a contact structure $\xi$ on $M$ is \textbf{supported by $\pi$} if $\xi$ is isotopic to a contact structure that admits a Giroux form for $\pi$.

\begin{remark}
    We will frequently use the $3$-tuple $(P, \phi, \Sigma)$ as a shorthand for a spinal open book where $P$, $\phi$, and $\Sigma$ denote the page, monodromy, and vertebra, respectively. Note that this tuple does not uniquely determine the spinal open book, as it omits the specific correspondence between vertebrae and the boundary components of the page.
\end{remark}

\begin{definition}\label{def:symmetricsobd}
    A spinal open book $\pi = (\pi_{\Sigma}, \pi_P)$ of $M$ is \textbf{symmetric} if
    \begin{enumerate}
        \item all pages are diffeomorphic,
        \item for each vertebra $\Sigma_i$ $(\Sigma = \Sigma_1 \sqcup \cdots \sqcup \Sigma_r)$, there exists $k_i \in \N$ such that every page has exactly $k_i$ boundary components in ${\pi_\Sigma}^{-1}(\partial \Sigma_i)$.
    \end{enumerate}
\end{definition}

We also say that $\pi$ is \textbf{uniform} if it is symmetric and there exists a compact oriented surface $B$ such that
\begin{enumerate}
    \item connected components of $\partial B$ are in bijection with the components of $M_P$;
    \item for each vertebra $\Sigma_i$ $(\Sigma = \Sigma_1 \sqcup \cdots \sqcup \Sigma_r)$, there exists a $k_i$-fold branched cover $\pi_i: \Sigma_i \to B$ such that for each boundary component $\gamma \subset \partial \Sigma_i$, $\pi_i|_{\gamma}$ is an $m_\gamma$-fold cover where $m_\gamma$ denotes the multiplicity of $\pi_P$ at $\pi_\Sigma^{-1}(\gamma) \subset M_P$.
\end{enumerate}

In the above case we further say that $\pi$ is \textbf{uniform with respect to $B$}, and we will refer to $B$ as the \textbf{base} of the spinal open book. We say that $\pi$ is \textbf{Lefschetz-amenable} if it is uniform and all branched covers have no branch points. These definitions are motivated by the fact that a bordered Lefschetz fibration induces a Lefschetz-amenable spinal open book on its boundary, while a nearly Lefschetz fibration induces a uniform spinal open book. See \S\ref{sec: PANLF} for the definition of nearly Lefschetz fibrations and see \S\ref{sec: examples} for examples of the construction of a nearly Lefschetz fibration from a uniform spinal open book.

The Lefschetz amenable case has been extensively studied in \cite{spinal_I,spinal_II}, while \cite{spinal_exotic} studies the non-Lefschetz-amenable case which necessarily involves the exotic fibers. See \S\ref{subsec: exotic} for more details on exotic fibers. In the non-amenable case, we will refer to a vertebra (and a spine) as {\bf non-trivial} if the map to $B$ has branch points.

We present some examples. For a Lefschetz-amenable spinal open book, see \S\ref{subsubsec: T_sobd} for a spinal open book for $3$-torus $\bT^3$. The next example describes a non-uniform spinal open book decomposition.

\begin{example}
    Let $P = P_1 \sqcup P_2$ be two sheets of pair of pants and $\Sigma = \Sigma_1 \sqcup \Sigma_2 \sqcup\Sigma_3$ where $\Sigma_1$ and $\Sigma_2$ are disks and $\Sigma_3$ is an annulus. The monodromy of the pages is $\phi = (\phi_1, \phi_2)$, where $\phi_1$ interchanges two inner boundary components of $P_1$ and $\phi_2$ interchanges two inner boundary components of $P_2$. See \S\ref{subsec: framed_mcg} for the definition of a boundary interchange map.
    To determine the spinal open book, we specify the correspondence between vertebrae and boundary components of the pages. First, $\Sigma_1$ and $\Sigma_2$ correspond to the outer boundary components of $P_1$ and $P_2$, respectively, and $\Sigma_3$ correspond to the inner boundary components of both $P_1$ and $P_2$.
    Here, the base $B$ should have two boundary components as there are two connected components in $P$. However, the disks $\Sigma_1$ and $\Sigma_2$ cannot be a (branched) cover of $B$, so the spinal open book cannot be uniform.
\end{example}

%-----------------------------------------------------------------------
\subsection{Framed mapping class group} \label{subsec: framed_mcg}
%-----------------------------------------------------------------------
To describe the monodromy of a spinal open book decomposition, we use the \textbf{framed mapping class group} of the page. In this section, we review relevant material on framed mapping class groups. For more details, readers are referred to \cite{baykur2016multisections}.

\begin{definition}
    Let $M$ be a manifold. Suppose $P \subset M$ and $U_1,\ldots,U_k \subset TM$ such that $\pi(\bigcup U_i) \cap P = \emptyset$.
    The \textbf{framed mapping class group} of $M$ relative to $P$ and $U_1,\ldots, U_k$ is defined as
    \begin{align*}
        \Mod_P(M; U_1, \ldots, U_k) &:= \pi_0(\Diff_P^+(M;U_1,\ldots,U_k))\\
        &:= \pi_0\left(\{\phi \in \Diff^+(M) \mid \phi|_{P} = \operatorname{id},\;d\phi(U_i)=U_i\}\right).
    \end{align*}
\end{definition}

For a compact oriented surface $\Sigma$, we use the following version of framed mapping class group:

\begin{definition}
    Let $\Sigma$ be a compact oriented surface with $n$ boundary components and $\{p_1,...,p_n\}$ a set of points one on each boundary component. In this case, we define its framed mapping class group as
    \[
        \FMod(\Sigma) := \Mod_{\emptyset}(\Sigma; \{p_1,\ldots,p_n\})
    \]
    Here, $U_1 := \{p_1,\dots, p_n\}$ and we identify each $p_i$ with the zero vector in $T_{p_i}\Sigma$.
\end{definition}

\begin{remark}
  The definition above was introduced by Baykur and Hayano \cite{baykur2016multisections}, which is equivalent to the mapping class group of a punctured surface that fixes the vectors at the punctures setwise. Thus, we can consider $\FMod(\Sigma)$ as the mapping class group that fixes the frames at the punctures.  %\ar{clarify this convoluted notion as a simple thing}
\end{remark}

\begin{remark}
    In the previous paper \cite{spinal_exotic}, we defined the \textbf{spinal mapping class group} $\operatorname{SMod}(\Sigma)$, which turned out to be identical to the framed mapping class group $\FMod(\Sigma)$. Consequently, we use the term ``framed mapping class group'', as it is a term already used in the literature.
\end{remark}

We now describe a set of generators for the framed mapping class group $\FMod(\Sigma)$. In particular, we note that the \textbf{half twist $\tau'_{\gamma}$} between two boundary components about an arc $\gamma$ is as in Figure~\ref{fig:monodromy1}. This map rotates and switches the two boundary components around a small neighborhood of the arc in a counterclockwise direction. For our purpose, however, we instead use a \textbf{boundary interchange} $\tau_{\gamma} = T_{c_1}^{-1} \circ T_{c_2}^{-1} \circ \tau'_{\gamma}$ that interchanges two boundary components and then negatively twists the boundary components to account for trivializations coming from spinal open book decompositions. See Figure~\ref{fig:boundarytwist}.

\begin{theorem}
    $\FMod(\Sigma)$ is generated by Dehn twists about homologically essential simple closed curves on $\Sigma$ and boundary interchange maps between two boundary components about an arc connecting two points on each boundary components.
\end{theorem}

%-----------------------------------------------------------------------
\subsection{Monodromy factorizations} \label{subsec: factorization}
%-----------------------------------------------------------------------
Here, we first introduce a monodromy factorization of a framed mapping class, and then explain how to apply it to a spinal open book setup. To do so, we recall from $\S\ref{subsec: framed_mcg}$ that a boundary interchange map $\tau_{\gamma}$ interchanges two boundary components of $P$ along an arc $\gamma$, as shown in Figure $\ref{fig:boundarytwist}$.

Let $P$ be a compact oriented surface with boundary, $\phi \in \FMod(P)$, and $\rho\colon G \to \Mod(P)$ a monodromy representation from a group $G$.

\begin{definition}
    Let $P, \phi,$ and $\rho$ be as above. We further suppose $g \in G$ and $n \in \mathbb{N}$. A \textbf{positive admissible factorization of $\phi$ with respect to $(\rho, g, n)$} is an expression of $\phi \circ \rho(g)^{-1}$ as a product of:
    \begin{itemize}
        \item $n$ boundary interchange maps between pairs of boundary components of $P$ along arcs that join the boundary components, and
        \item positive Dehn twists about homologically essential curves on $P$.
    \end{itemize}
\end{definition}

See Figure $\ref{fig:ellipticbundles2}$ for an example of such a factorization.

Suppose we have a spinal open book $(\pi_{\Sigma},\pi_P)$ of a closed $3$-manifold $M$ with a page $P = (P_1, \ldots, P_m)$, monodromy $\phi = (\phi_1,\ldots, \phi_m)$ where $\phi_i \in \FMod(P)$, and vertebrae $\Sigma = \Sigma_1 \sqcup \dots \sqcup \Sigma_r$. Let $B$ be the base of the spinal open book, a compact oriented surface with boundary such that for each $k = 1,\ldots,r$, there is a (possibly branched) covering map $\pi_k :\Sigma_k \to B$ with $n_k$ simple branch points. Also suppose $(\pi_{\Sigma},\pi_P)$ is uniform with respect to $B$. Then all pages are diffeomorphic and there is a one-to-one correspondence between each paper component and boundary component of $B$. To be more specific, for each boundary component $\partial_i$ of $B$, there exists a corresponding paper component with monodromy $\phi_i$ that intersects the boundary component $\pi_k^{-1}(\partial_i)$ of $\Sigma_k$ for every $k=1,\ldots,r$.

Now we consider the following subsets of $\partial P$:
\begin{align*}
    \partial_k P_i &:= \{ c \subset \partial P_i \mid c \text{ is a connected component intersecting a spine component } \Sigma_k \times S^1 \}
\end{align*}

From now on, we will slightly abuse notation. We first set $P_0$ as a ``model surface'' such that there is a diffeomorphism $f_i \colon P_0 \to P_i$ for each $i = 1,\ldots,m$. Then we can identify all pages $P_1,\ldots,P_m$ with the single surface $P_0$ and use the monodromy $f_i^{-1} \circ \phi_i \circ f_i$ instead of $\phi_i$. We relabel $f_i^{-1} \circ \phi_i \circ f_i$  as $\phi_i$ to shorten the notation. Also, we fix a monodromy representation $\rho\colon \pi_1(B) \to \Mod(P_0)$. Recall that for each paper component with monodromy $\phi_i$, there is a corresponding boundary component $\partial_i$ of $B$.

\begin{definition}\label{def: pos_adm_fact}
    Let $P$, $\phi$, $B$, $\rho$, and $n_k$ be as above. A {\bf positive admissible factorization of $\phi$ with respect to $(B, \rho)$} is an expression of each $\phi_i \circ {\rho}(\partial_i)^{-1}$ for $i=1,\ldots,m$ as a product of
    \begin{itemize}
        \item $n_k$ boundary interchanges between pairs of boundary components in $\partial_k P_i$ along arcs that join the boundary components for $k = 1,\dots,r$, and
        \item positive Dehn twists about homologically essential curves in $P_i$.
    \end{itemize}
\end{definition}

The main purpose of this factorization is to construct a PANLF over $B$ with fiber $P_0$. We first construct a $P$-fiber bundle over $B$ using the monodromy representation $\rho$. We then introduce exotic fibers by adding Weinstein 2-handles corresponding to boundary interchanges. Lastly, add singular fibers by attaching Weinstein $2$-handles along the curves associated with the Dehn twists. Note that each connected component of paper $M_P$, a $P_i$-bundle over $S^1$ with monodromy $\phi_i$, corresponds to a subbundle of this PANLF restricted to the boundary component $\partial_i$ of $B$. The total number of exotic fibers is $n_B = n_1 + \cdots + n_r$, and the total number of singular fibers is equal to the total number of Dehn twists across $P_1,\ldots,P_m$.

Now we collect all such factorizations that range over all monodromy representations.

\begin{definition}
    The collection of positive admissible factorizations of $\phi$ with respect to $(B, \rho)$, as $\rho$ ranges over all possible monodromy representations, will be called the {\bf set of  positive admissible factorizations of $\phi$ with respect to $B$}.
\end{definition}

\begin{example}
    Figure~\ref{fig:ellipticbundles} represents spinal open book decompositions for elliptic torus bundles. However, the factorizations of the monodromies $\phi_1$, $\phi_2$, and $\phi_3$ are not positively admissible, as they contain negative Dehn twists along closed curves. Figure~\ref{fig:ellipticbundles2} represents the same elliptic torus bundles, but their monodromy factorizations are positively admissible.
\end{example}

We now describe equivalence between two factorizations that induce an equivalent PANLF. The following definition is due to Baykur and Hayano \cite{Baykur_Hayano_Hurwitz}.

\begin{definition}[Generalized Hurwitz equivalence] \label{def: generalized_hurwitz_equivalence}
    Let $F$ be a compact oriented surface with boundary. Two positive allowable monodromy factorizations of $\phi \in \FMod(F)$ will be said to be \textbf{Hurwitz equivalent} if they are related by the following moves:
    \begin{enumerate}
        \item {\bf Elementary transformation}, which changes a factorization as follows:
        \[
            \eta_{k+l}\cdots\eta_{i+1}\eta_i\cdots\eta_1 \longleftrightarrow \eta_{k+l}\cdots(\eta_{i+1}\eta_i\eta_{i+1}^{-1})\eta_{i+1}\cdots\eta_1
         \]
        \item {\bf Global conjugation}, which changes each member of a factorization by the conjugation of some mapping class $\psi \in \Mod_\partial(F)$:
        \[
            \eta_{k+l}\cdots\eta_1 \longleftrightarrow (\psi\eta_{k+l}\psi^{-1})\cdots(\psi\eta_1\psi^{-1})
        \]
        \item {\bf Framing conjugation}, which changes a factorization as follows:
        \[
            \eta_{k+l}\cdots\eta_{i+1}\eta_i\eta_{i-1}\cdots\eta_1 \longleftrightarrow \eta_{k+l}\cdots\eta_{i+1}(T_{\delta}\eta_iT_{\delta}^{-1})\eta_{i-1}\cdots\eta_1
        \]
        where $\delta$ is a simple closed curve parallel to a boundary component of $F$.
    \end{enumerate}
\end{definition}

The generalized Hurwitz moves naturally arise as a consequence of choices made when obtaining a factorization from a PANLF as we will see in \S \ref{sec: PANLF}. As a result, we obtain the following result:

\begin{proposition}
    Two Hurwitz equivalent monodromy factorizations give rise to two symplectic PANLFs that are deformation equivalent.
\end{proposition}

%-----------------------------------------------------
\subsection{Construction of spinal open book decompositions}\label{subsec:spinal_obd_examples}
%-----------------------------------------------------
In this section we discuss some methods to construct supporting spinal open book decompositions for contact 3-manifolds.

%-----------------------------------------------------
\subsubsection{The Stein fillable contact structure on the 3-torus} \label{subsubsec: T_sobd}
\begin{figure}
    \centering
    \includegraphics[width=0.7\linewidth]{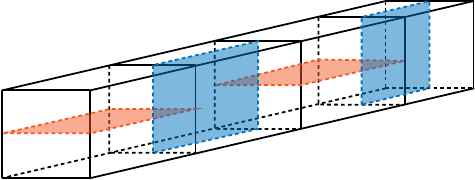}
    \caption{The 3-torus $\bT^3$. The top and bottom faces are identified, as are the front and back, and the left and right. The blue annuli represent the pages of the spinal open book, and the red annuli represent the vertebrae.}
    \label{fig:T3}
\end{figure}

%-----------------------------------------------------
Eliashberg \cite{Eliashberg:T3} showed that there exists a unique Stein fillable contact structure $\xi_0$ on $\bT^3$. In \cite{Wendl_hierarchy, spinal_II}, a spinal open book decomposition for $(\bT^3, \xi_0)$ was introduced, as shown in Figure~\ref{fig:T3}. The first and third $\bT^2 \times I$ pieces in the figure are spines and the second and fourth are papers. The blue annuli are pages, and the red annuli are vertebrae. This spinal open book decomposition naturally arises as the boundary of the bordered Lefschetz fibration $A \times A$, with pages $\partial A \times A$ and vertebrae $A \times \partial A$, where $A = S^1 \times I$ is an annulus. The monodromy of the paper is $\operatorname{id}$.

Now we verify that this spinal open book decomposition actually supports $\xi_0$. According to \S\ref{subsec: bundles}, a spine with an annulus vertebra is equivalent to a relative open book deposition with the word $(aba)^{-2}$ as shown in the shaded region in Figure~\ref{fig:sobd_parabolic}. Since our spinal open book consists of two annulus pages with identity monodromy, and two annulus vertebrae, it is equivalent to an ordinary open book decomposition with the word $(aba)^{-4}$, which supports $\xi_0$ by Van Horn-Morris \cite{Van_horn_morris_thesis}.

%--------------------------------------------------
\subsubsection{Blown up summed open books}
%--------------------------------------------------
There is a natural operation on a regular contact open book decomposition called \textbf{binding sum} that gives rise to a spinal open book. Let $\pi: M\setminus B\to S^1$ be an open book decomposition with at least two binding components $B_1, B_2 \subset B$. Topologically, one removes neighborhoods of $B_1$ and $B_2$ and glues the resulting boundary tori via an orientation reversing diffeomorphism which maps the boundaries of the pages to each other and meridians to meridians (with opposite orientation). This is equivalent to attach two boundary components to each boundary component of $\bT^2 \times I$. Thus the resulting manifold naturally has the structure of a spinal open book decomposition with an annulus vertebra replacing the two disk vertebrae, and the supporting contact structure is obtained by a contact fiber sum along $B_1$ and $B_2$ of the original open book decomposition.

%--------------------------------------------------
\subsubsection{Torus bundles and relative open books} \label{subsec: bundles}
%--------------------------------------------------
Here, we construct planar spinal open book decompositions for \textbf{rotational contact structures} on torus bundles, using the technique of \textbf{relative open book decompositions} introduced by Van Horn-Morris \cite{Van_horn_morris_thesis}.

we first review rotational contact structures on torus bundles. Let $A \in SL_2(\Z) = \Mod(\bT^2)$. Then we define the torus bundle with monodromy $A$ by
\begin{align*}
    T_A := \bT^2 \times \R / (x,t) \sim (Ax,t-1).
\end{align*}
Moreover, $T_A$ is called {\bf elliptic, parabolic}, and {\bf hyperbolic}, respectively, as $|\mathrm{tr}(A)|$ is less than 2, equal to 2, or greater than 2. We consider a {\bf rotational contact structure} $\xi_n$ for an integer $n \geq 0$ defined by
\begin{align*}
    \xi_n = \ker[\sin \phi_e(t)\,dx + \cos \phi_n(t)\,dy],\;\; (x,y) \in \bT^2, \,t \in \R / \Z
\end{align*}
where $\phi_n \colon \R/\Z \to \R$ satisfies
\[
    n\pi \leq \sup_{t \in \R / \Z}(\phi_n(t+1) - \phi_n(t)) < (n+1)\pi.
\]
We say $\xi_n$ has \textbf{$n \pi$-twisting}. Each torus bundle only admits either odd or even rotational contact structures. If $A$ is parabolic, then $A$ is conjugate to
\[
    \pm A_k = \pm\begin{pmatrix}1 & k \\ 0 & 1\end{pmatrix}
\]
for some $k \in \Z$. We say a parabolic torus bundle $T_A$ is {\bf positive} if $A$ is conjugate to $A_k$ and {\bf negative} if $A$ is conjugate to $-A_k$. Let $T_+(k) := T_{A_k}$ and $T_-(k) := T_{-A_k}$. Note that $T_{+}(k)$ only admits even rotational contact structures $\xi_{2n}$ for $n \geq 1$ and $T_-(k)$ only admits odd rotational contact structures $\xi_{2n+1}$ for $n \geq 0$.

Next, we review relative open book decompositions, introduced by Van Horn-Morris \cite{Van_horn_morris_thesis}. Let $P$ be a compact oriented surface with boundary and $\phi \in \Mod_{\partial}(P)$. Consider the mapping torus
\[
    P_{\phi} = P \times [0,1] / \sim
\]
where $(p,1) \sim (\phi(p),0)$ for $p \in P$. Now pick a partition of $\partial P$ into two sets $\partial_B$ and $\partial_T$. Each component of $\partial_B$ is called {\bf binding circles}, and each component of $\partial_T$ is called {\bf boundary circles}. We foliate the boundary tori $\partial_B \times S^1 \subset \partial P_{\phi}$ into $\{x\} \times S^1$, where $x \in \partial_B$. From this, we can construct a $3$-manifold $M$ with torus boundary components by collapsing the foliating circles in each component of $\partial_B \times S^1$. We say $(P,\phi,(\partial_B,\partial_T))$ is a {\bf relative open book decomposition} of $M$. If $\partial_T = \varnothing$, we obtain a closed manifold with a regular open book decomposition.

As in the case of regular open book decompositions, there is a contact structure compatible with a given relative open book decomposition.

\begin{definition}
    Let $M$ be a compact oriented $3$-manifold and $(P,\phi,(\partial_B,\partial_T))$ a relative open book of $M$. A contact structure $\xi$ on $M$ is said to be {\bf compatible with} $(P,\phi,(\partial_B,\partial_T))$ if there exists a contact form $\alpha$ where $d\alpha$ restricts to a symplectic form on each page, and the foliation given by $\xi$ on $\partial M$ agrees with the foliation given by boundary circles in $\partial_T$.
\end{definition}

\begin{figure}[htbp]{
    \vspace{0.2cm}
    \begin{overpic}[tics=20]
    {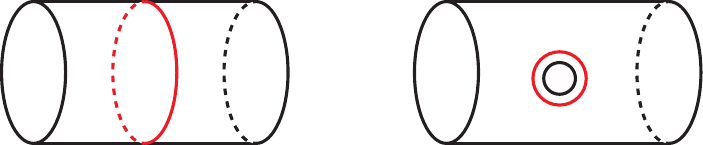}
    \end{overpic}}
    \vspace{0.2cm}
    \caption{Relative open book decompositions for $a$ and $b^{-1}$. For $a^{-1}$ and $b$ we use the negative Dehn twists instead.}
    \label{fig:relativeBook}
\end{figure}

We now review a method for representing relative open books of $S^1 \times I \times S^1$ using \textbf{words}. Four relative open books, which correspond to the set of letters $\{a,a^{-1},b,b^{-1}\}$ as depicted in Figure~\ref{fig:relativeBook}, act as building blocks for a relative open book decomposition supporting a contact structure on $\bT^2 \times I$, or a torus bundle over $S^1$. A torus bundle obtained by gluing these pieces is denoted by the corresponding word in terms of $a$, $a^{-1}$, $b$ and $b^{-1}$. For example, the left drawing of Figure~\ref{fig:sobd_parabolic} represents an open book decomposition corresponding to the word $a^k(aba)^{-4}$.

Lastly, we find a relative open book of a spine component with annular vertebra that agrees on the boundary of the spine component with the page foliation of the spinal open book, which is a key to construct planar spinal open book decompositions for rotational torus bundles. Note that a spine component with an annulus vertebra is also homeomorphic to $S^1 \times I \times S^1$.

\begin{proposition}[Min-Roy-Wang \cite{spinal_exotic}] \label{prop:obd_ann_vertebra}
    The relative open book for a spine component of an annular vertebra is contactomorphic to a relative open book with the word $(aba)^{-2}$.
\end{proposition}

Using this proposition, we can convert any genus-$1$ open book decomposition containing a word $(aba)^{-2}$ into a planar spinal open book decompositions as follows:

\begin{enumerate}[(i)]
    \item Consider a genus-$1$ open book decomposition that contains the word $(aba)^{-2}$.

    \item Remove the relative open book corresponding to $(aba)^{-2}$ from the entire open book decomposition. Then the remaining part become a page of the spinal open book, and the relative open book becomes an annular vertebra.

    \item Repeat (ii) on the page of the spinal open book and we obtain multiple annular vertebrae and page components. This new spinal open book supports the same contact structure as the original open book.
\end{enumerate}

See Figure~\ref{fig:sobd_parabolic} for an example. An open book decomposition with the word $(aba)^{-4} a^k$ supports a positive parabolic torus bundle $(T_+(k),\xi_2)$. Now we remove two relative open books corresponding to $(aba)^{-2}$ and replace them with annular vertebrae. The remaning part becomes two annular pages with monodromies $a^k$ and $\operatorname{id}$, respectively.

\begin{exercise} \label{exer:elliptic}
    Find planar spinal open book decompositions for $\pi$-twisting elliptic torus bundles.\\
    ({\em Hint: Find the words for the bundles from \cite[Theorem~4.3.1]{Van_horn_morris_thesis}}.)
\end{exercise}

\begin{exercise} \label{exer:3-torus}
    Show that the $2\pi$-twisting $3$-torus has the word $(aba)^{-4}$.\\ {\em (Hint: Convert a spinal open book for the $3$-torus from \S\ref{subsubsec: T_sobd} into an honest open book.)}
\end{exercise}

\begin{exercise}
    Show that any rotational torus bundle with at least $\pi$-twisting admits a planar spinal open book decomposition. \\{\em (Hint: Use Exercise~\ref{exer:3-torus} to show that these bundles admit a genus-1 open book containing the word $(aba)^{-2}$.)}
\end{exercise}

%%%%%%%%%%%%%%%%%%%%%%%%%%%%%%%%%%%%%%%%%%%%%%%%
\section{Nearly Lefschetz fibrations}\label{sec: PANLF}
%%%%%%%%%%%%%%%%%%%%%%%%%%%%%%%%%%%%%%%%%%%%%%%%

In this section we describe the topological structure on a symplectic 4-manifold with boundary that induces a spinal open book on the boundary. Just as bordered Lefschetz fibrations over the disk induce open book decompositions on the boundary, these \textbf{nearly Lefschetz fibrations} induce spinal open book decompositions of the boundary. We describe the notion of a bordered Lefschetz fibration over a surface with boundary in \textsection~\ref{subsec: blf}. These give rise to ``amenable'' spinal open books on the boundary. Then, in \textsection~\ref{subsec: exotic}, we describe the ``singularity-at-infinity'' that arises when studying the foliations coming from planar spinal open books, called \textbf{exotic fibers}. We give an explicit local model and study the induced monodromy on an $S^1$-bundle over a curve in the base enclosing such a these singularity. Then in \textsection~\ref{subsec: panlf_stein}, we describe nearly Lefschetz fibrations, and how to associate Stein structures to them.

%-----------------------------------------------------
\subsection{Bordered Lefschetz fibrations}\label{subsec: blf}
%-----------------------------------------------------
A spinal open book decomposition naturally arises as the boundary of a bordered Lefschetz fibration. In this section we recall the definition of a bordered Lefschetz fibration by following the conventions of \cite{spinal_I}. Let $E$ be a compact oriented connected 4-manifold with corners. We decompose the boundary of $E$ into two parts along the corners:
\[
    \partial E = \partial_h E \cup \partial_v E.
\]
We further assume that the intersection of $\partial_h E$ and $\partial_v E$ is a collection of $2$-tori.

Unlike \cite{spinal_I}, since we will not be carrying out any analysis near the corners, we will treat $E$ as a smooth manifold, assuming we already rounded the corners.
Let $B$ denote a compact oriented connected surface with connected boundary.

\begin{definition}
    A \textbf{bordered Lefschetz fibration} of $E$ over $B$ is a smooth map $\Pi: E \to B$ with finitely many interior critical points $E^{\crit} \subset E^{\circ}$ and critical values $B^{\crit} \subset B^{\circ}$ such that
    \begin{enumerate}
        \item  $\Pi^{-1}(\partial B) = \partial_v E$ and $\Pi|_{\partial_v E} :\partial_v E \to \partial B$ is a smooth fiber bundle. %This corresponds to the paper of the boundary spinal open book.
        \item $\Pi|_{\partial_h E} :\partial_h E \to B$ is a smooth fiber bundle. %This corresponds to the spine on the boundary spinal open book.
        \item $(\Pi|_{\partial_h E}, \Pi|_{\partial_v E})$ is a spinal open book decomposition of $\bd E$.
        \item For each $p \in E^{\crit}$ and $\Pi(p) \in B^{\crit}$, there are local holomorphic coordinates on $E$ centered at $p$ and on $B$ centered at $\Pi(p)$ such that $\Pi(z_1, z_2) = z_1^2 + z_2^2$.
        \item For each $z\in B$, the fiber $E_z := \Pi^{-1}(z)$ is connected and has nonempty boundary in $\partial_h E$.
    \end{enumerate}
\end{definition}

The fibers $\Pi^{-1}(z)$ for $z \in B \setminus B^{\crit}$ are called \textbf{regular fibers}, and for $z \in B^{\crit}$ are called \textbf{singular fibers}. The regular fibers are all homeomorphic smooth compact oriented surfaces with boundary. The singular fibers are smoothly immersed connected surfaces with positive transverse self-intersections -- topologically they are obtained by taking a regular fiber and pinching a curve on it down to a point; the curve that gets pinched is called a \textbf{vanishing cycle}. A Lefschetz fibration $\Pi$ is \textbf{allowable} if all the vanishing cycles are homologically essential curves on a regular fiber.

It is well-known that allowable bordered Lefschetz fibrations admit a Stein structure.

\begin{theorem}[{\cite[Theorem 3.9]{baykurvhm}\cite[Theorem B]{spinal_I}}]\label{thm: blf_stein}
    Let $E$ be a 4-dimensional allowable bordered Lefschetz fibration and $\pi$ is an induced spinal open book of $\partial E$. Then, $E$ admits a (well-defined up to Stein homotopy) Stein structure which is a Stein filling of the contact manifold $(\partial E, \pi)$. Moreover, given two bordered Lefschetz fibrations $E_1$ and $E_2$ that fill $(\partial E, \pi)$, their Stein structures can be chosen to induce the same contact structure on the boundary.
\end{theorem}

%-------------------------------
\subsection{Exotic fibers}\label{subsec: exotic}
%-------------------------------

In this section, we will recall the local model of a neighborhood of an exotic fiber and explain how to see the local monodromy in terms of a boundary interchange in Lemma \ref{lem:monodromy}. This should be reminiscent of the story Lefschetz singularities, where we have holomorphic local models $(z_1, z_2) \mapsto z_1^2+z_2^2$ and the local monodromy is a right-handed Dehn twist around an homologically essential curve. See \cite{gompf_stipsicz}. However, it should be emphasized early on that, as opposed to the Lefschetz singular fibers, the exotic fibers exhibit a ``singularity-at-infinity'' phenomenon, incorporating a mapping class that interchanges boundary components.

The original definition of exotic fibers is a technical one: they are the main level of buildings which appear in the compactification of certain moduli space of curves that \cite{spinal_II} analyzes. Since this lecture note focuses on the application of spinal open books and symplectic fillings potentially containing exotic fibers, we will focus on the topological description and refer the readers to \cite[\textsection 4]{spinal_exotic} for a careful treatment of the local model of the exotic fibers in terms of SFT moduli spaces and the proof of uniqueness of the local model.

\begin{definition}[Local model] \label{def:model}
    Consider the sequence of curves
    \[
        u_c: \mathbb{C}\backslash\{\pm 1\} \to (\mathbb{C}\backslash\{0\}) \times \mathbb{C},
    \]
    defined by
    \begin{equation} \label{eq:local_model}
        u_c(z)= (c(z^2-1), \sqrt{c}z).
    \end{equation}
    Consider the projection map
    \begin{equation}\label{eq:lefschetz_projection}
        \Pi_0: (\mathbb{C}\backslash\{0\})\times \mathbb{C} \to \mathbb{C} \text{ sending } (z_1,z_2)\mapsto z_2^2-z_1,
    \end{equation}
    whose level sets are identified with the image of $u_c$.
\end{definition}

\begin{figure}[htbp]{
  \vspace{0.2cm}
  \begin{overpic}[tics=20]
  {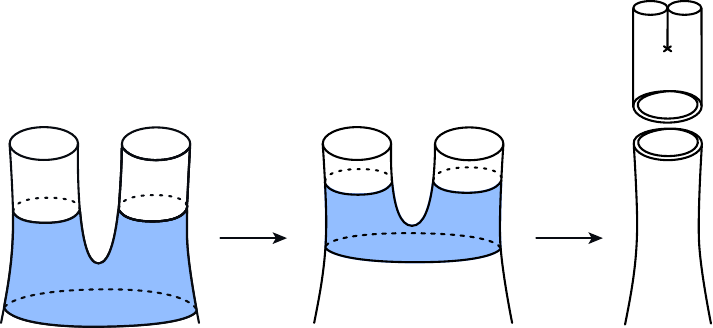}
  \end{overpic}}
  \vspace{0.2cm}
  \caption{A local description of the SFT limit from regular fibers to an exotic fiber.}
  \label{fig:building}
\end{figure}

\begin{remark}
    The moduli space of curves $u_c(z)$ in Equation~(\ref{eq:local_model}) is parameterized by the domain of $c\in \C\setminus\{0\}$. The reparameterization by $z':= \sqrt{c}z$ allows us to compactify the moduli space at $c = 0$ with a limiting building as in Figure \ref{fig:building}.
\end{remark}

\begin{remark} \label{rmk: compactified_model}
    Definition \ref{eq:local_model} as well as the understanding of pseudoholomorphic foliations in \cite{spinal_II} are both in the non-compact setting in the completion of the minimal strong symplectic fillings. As is often useful to talk about symplectic fillings as a compact manifold with contact boundary, the corresponding local model also needs to be compactified to be technically correct. We define the \textbf{compactified local model} to be the restriction of the above model to $(\C \setminus B) \times \C$, which corresponds to cutting off the cylindrical ends at the punctures. Here $B$ is defined to be the closed $\epsilon$-disk about the origin in $\C$. The restricted foliation on $(\C \setminus B) \times \C$ is as follows. For $u_c$, $\im(u_c) \cap (\C \setminus B) \times \C = \{u_c(z) | |z^2-1| > \frac{\epsilon}{|c|}\}$. Thus, for $|c| > \epsilon$,
    $\im(u_c) \cap (\C \setminus B) \times \C$ is topologically a pair of pants, while for $|c| < \epsilon$ $\im(u_c) \cap (\C \setminus B) \times \C$ is topologically an annulus, with the curves corresponding to $|c| = \epsilon$ having a singularity. So in the compactified local model, we get \textbf{regular curves}, corresponding to $\im(u_c)$ for $c > \epsilon$, and a $\mathbb{D}^2$-worth of \textbf{exotic curves}, corresponding to $\im(u_c)$ for $c \leq \epsilon$, which we call an \textbf{exotic neighborhood}.
\end{remark}

\begin{figure}[htbp]{
  \vspace{0.2cm}
  \begin{overpic}[tics=20]
  {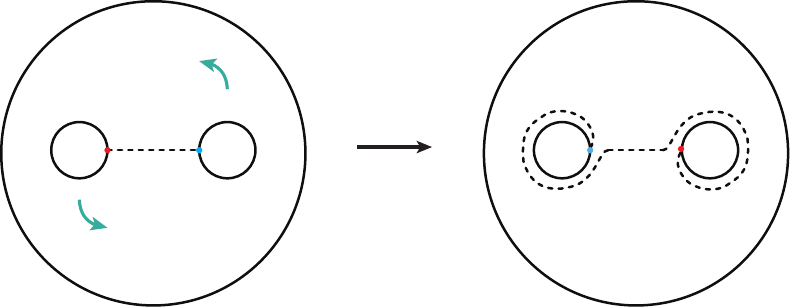}
   \put(36,95){\scriptsize$1$}
   \put(108,95){\scriptsize$2$}
   \put(74,85){\scriptsize$\gamma$ }

   \put(269, 98){\scriptsize$2$}
   \put(342, 98){\scriptsize$1$}
   \put(306, 85){\scriptsize$\gamma$}
  \end{overpic}}
  \vspace{0.2cm}
  \caption{A local description of a \textbf{boundary interchange} $\tau_{\gamma}$, where $\gamma$ is the equatorial arc connecting the two boundary components in the interior.}
  \label{fig:boundarytwist}
\end{figure}

We will see in Lemma~\ref{lem:monodromy} that this gives a model that interchanges the two punctures as one goes around $0\in \mathbb{C}$ as in Figure \ref{fig:boundarytwist}. Further, this local model of the exotic fiber also allows us to show that the number of exotic points is equal to the number of the branch points of $\Pi|_{\Sigma}$ in Theorem \ref{thm: intro-localmodel}. To describe the monodromy of the fibers in a compactified filling around the exotic neighborhood, we will talk about \textbf{boundary components} of the fibers in the compactified model, which correspond to punctures in the local model.

\begin{lemma}\label{lem:monodromy}
    Consider the fibers in the compactified local model. Outside the exotic neighborhood, the fibers are pairs of pants. The monodromy of the bundle with pairs of pants fibers over $S^1$, around the exotic fiber $u_0$, is given by interchanging two boundary components counter-clockwise along an arc and simultaneously rotating the boundary components by $\pi$ clockwise.
\end{lemma}

\begin{proof}
    In the model from Definition~\ref{def:model}, consider the following disk neighborhood $V$ of $0\in \bC$. Let $\delta>2\epsilon$ be large enough so that the $\delta$-disk neighborhood $ V\subset \bC$ contains the image of the exotic neighborhood. We parameterize $V$ so that $\Pi^{-1}(0)$ corresponds to the exotic fiber. Let $\zeta:= \partial V$. Over each point $x\in\zeta$, $\Pi^{-1}(x)$ is a pair of pants. We will characterize the monodromy of this bundle with pair of pants fibers around $\zeta$. Also consider a chart $U$ of the filling $W$ around the exotic neighborhood, as above. By construction, topologically $\Pi^{-1}(V) - U$ is a trivial $\Sigma - P$ bundle over the disk $V$, where $P$ is a pair of pants. In other words, the monodromy is supported in $P$.

    Now we parameterize $\zeta = \{\delta e^{2\pi i t}| t \in [0,1]\}$ and understand the monodromy around it. First, fix a trivialization on the fibers $\Pi^{-1}(\delta e^{2\pi i t})$, such that the bundle over $\zeta$ is trivial outside $P$. We can parameterize the fibers in $\Pi^{-1}(\zeta) \cap U$ as in Definition~\ref{def:model}, i.e., the trivialization is given by
    \begin{equation}
        z \mapsto (e^{2\pi i t}(z^2-1),e^{\pi i t}z)
    \end{equation}
    for $z \in \Pi^{-1}(\delta) \cap U$. Projecting to the second coordinate, the trivialization we need to consider is $H_t$ such that
    \[
        H_t: z \mapsto e^{\pi i t}z.
    \]
    We can parameterize $P$ as $\{|z|^2\leq2\} \cap \Pi^{-1}(\delta) \cap U$, so outside that we need to modify the above trivialization so that the bundle over $\zeta$ is trivial. Consider the function $f$ as shown in Figure~\ref{fig:trivialization}, and consider the trivialization
    \begin{equation}\label{eqn:monodromy_trivialization}
        h_t: z \mapsto e^{f(|z|)\pi i t}z.
    \end{equation}

    \begin{figure}[htbp]{
        \vspace{0.2cm}
        \begin{overpic}[tics=20]
        {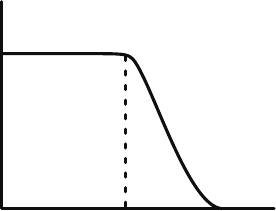}
        \put(-9, 74){\scriptsize$1$}
        \put(54, -10){\scriptsize$1.5$}
        \put(100, -10){\scriptsize$2$}
        \put(82, 100){\scriptsize$f(z)$}
        \end{overpic}}
        \vspace{0.2cm}
        \caption{The bump function needed to trivialize the monodromy for large $|z|$.}
        \label{fig:trivialization}
    \end{figure}

    \begin{figure}[htbp]{
        \vspace{0.2cm}
        \begin{overpic}[tics=20]
        {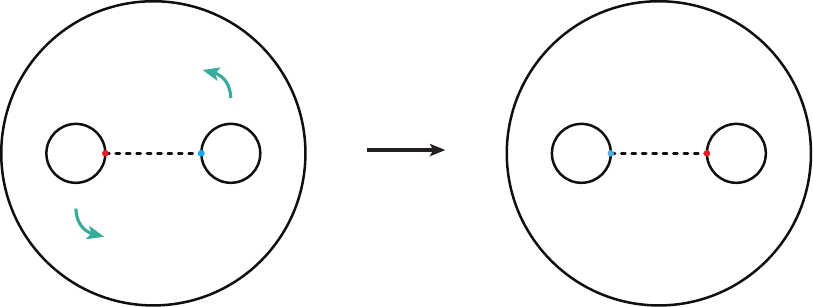}
        \put(35, 92){\scriptsize$1$}
        \put(109, 92){\scriptsize$2$}
        \put(72, 83){\scriptsize$\gamma$}

        \put(279, 92){\scriptsize$2$}
        \put(353, 92){\scriptsize$1$}
        \put(317, 83){\scriptsize$\gamma$}
        \end{overpic}}
        \vspace{0.2cm}
        \caption{A \dfn{half twist along $\gamma$}. Comparing the trivializations $h_0$ and $h_1$ after making the bundle trivial away from the pair of pants neighborhood.}
        \label{fig:monodromy1}
    \end{figure}

    Now, compare the trivializations $h_1$ and $h_0$. This is shown in Figure~\ref{fig:monodromy1} for the compactified model.

    The above modification of the trivializations ensures that we have localized the monodromy into the pair of pants neighborhood on a fiber. But we are not done yet, since to understand the monodromy of the bundle, in the sense of (3) and (4) in Definition~\ref{def:spinalobd}, we further need to trivialize the bundle near the boundaries of the fibers, to agree with the model of $\pi_P$ near this component of $\partial M_P \subset \Pi^{-1}(\zeta)$. Observe that the component of $\partial M_P$ in this local model has multiplicity 2, again in the sense of Definition~\ref{def:spinalobd}. Let us refer to the component of $\partial M_P$ as $T$, and give it coordinates $(\phi, \theta)$, where $\theta$ is the $S^1$-coordinate for the boundary component on a page. Locally $\pi_P$ has the form $\pi_P|_T(\phi, \theta) = 2\phi$. However, the trivializations in Equation~\ref{eqn:monodromy_trivialization} do not agree on $T$, for different fibers and their intersections with $T$. To make them agree, a clockwise $2\pi$-rotation of each of the boundary components of $P$ needs to happen. Using a bump function as before, which is 0 on the boundary components and 1 just outside, the trivialization $h_1$ looks like as in Figure~\ref{fig:boundarytwist}. By an isotopy, we can see that the monodromy is exactly a positive half-twist and a half clockwise twist of each of the boundary components, as described in Figure~\ref{fig:boundarytwist}.
\end{proof}

In \cite{spinal_exotic}, we show that this local model is unique up to symplectic deformation. This involves understanding the asymptotic neighborhood of finite-energy pseudoholomorphic curves asymptotic to a doubly-covered elliptic orbit. Theorem \ref{thm: main} is then a direct result of this local model, together with the curve classification in \cite[Proposition 1.30]{ spinal_II} in terms of regular, singular, and exotic fibers. In addition, the local model gives us the following.

\begin{mainthm}\label{thm: intro-localmodel}
    Let $(W, \omega)$ be a minimal strong filling of a contact manifold $(M,\xi)$ supported by a planar spinal open book. Then, the number of exotic fibers is equal to the number of branch points of $\Pi|_{S_i}$, where $\Pi$ is the same Lefschetz fibration as in Theorem \ref{thm: main}. Furthermore, the monodromy around an exotic fiber in a compactified filling is a boundary interchange on the nearby regular fibers, given by Figure~\ref{fig:boundarytwist}.
\end{mainthm}

This implies that $(W,\omega)$ has the structure of a nearly Lefschetz fibration, whose definition we will see now in the following subsection.

%----------------------------------------------------------
\subsection{PANLF and Stein structures}\label{subsec: panlf_stein}
%----------------------------------------------------------
For a symplectic filling of an honest planar open book decomposition with suitable conditions, a pseudoholomorphic foliation gives a Lefschetz fibration \cite{Wendl_strongly}. However, symplectic fillings of a general spinal open book admit a generalized version of a Lefschetz fibration, which we call a \textbf{nearly Lefschetz fibration}. In particular, it gives rise to a uniform spinal open book on its boundary, which is non-Lefschetz-amenable in general.

\begin{definition}\label{def: panlf}
  A \textbf{nearly Lefschetz fibration} of $E$ over $B$ is a smooth map $\Pi: E \to B$ with finitely many interior critical points $E^{\crit} \subset E^{\circ}$ and critical values $B^{\crit} \subset B^{\circ}$, and finitely many \textbf{exotic points} $B^{\exot} \subset B^{\circ}$ such that:
  \begin{enumerate}
    \item $\Pi^{-1}(\partial B) = \partial_v E$ and $\Pi|_{\partial_v E} :\partial_v E \to \partial B$ is a smooth fiber bundle.
    \item For each $p \in E^{\crit}$ and $\Pi(p) \in B^{\crit}$, there are local holomorphic coordinates on $E$ centered at $p$ and on $B$ centered at $\Pi(p)$ such that $\Pi(z_1, z_2) = z_1^2 + z_2^2$.
    \item All fibers $E_z := \Pi^{-1}(z)$ for $z\in B$ are connected and have nonempty boundary in $\partial_h E$.
    \item For each $c \in B^{\exot}$, there are local holomorphic coordinates $\psi \colon U_c \to \C$, where $U_c \subset B^{\circ}$, sending $c$ to $0$ and $\phi\colon V_c \to (\C - D) \times \C$, where $V_c \subset N(\partial_h (E)) \cap \Pi^{-1}(U_c)$ and $D$ is a small neighborhood of $0 \in \C$, such that the following diagram commutes
    \[
    \begin{tikzcd}
      V_c \arrow{r}{\phi}  \arrow{d}{\Pi} &(\C - D) \times \C\arrow{d}{\Pi_c}\\
      U_c \arrow{r}{\psi} & \C
    \end{tikzcd}
    \]
    where $\Pi_c(z_1,z_2) = z_2^2-z_1$.
    \item The boundary $\partial E = M$ admits a spinal open book decomposition, such that $\partial_vE = M_P$, $\partial_hE = M_\Sigma$, $\pi_P = \Pi|_{\partial_v E}$, and (unions of\footnote{Recall that by the definition of spinal open books, the fibers of the map $\pi_\Sigma$ has to be connected. Therefore, here the disjoint union of the fibers of $\pi_\Sigma$ corresponds to (possibly disconnected) fibers of $\Pi|_{\partial_h E}$.}) fibers of $\pi_\Sigma$ are the fibers of $\Pi|_{\partial_h E}$, away from the $V_c$ corresponding to $c \in B^{\exot}$.
  \end{enumerate}
\end{definition}

In the case of a nearly Lefschetz fibration, we have \textbf{regular} and \textbf{singular fibers} as in the case of a bordered Lefschetz fibration. The difference lies in the exotic points. We understood the local model of exotic points in the previous subsection. To reiterate, let $c\in B^{\exot}$ and $z \in U_c \subset B^o$, we have that $\Pi^{-1}(z)$ has the topology of a regular fiber near the boundary of $U_c$, but has one fewer boundary component in a neighborhood of the origin in $U_c$. We call all the fibers that have one fewer boundary component the \textbf{exotic fibers}\footnote{In \cite{spinal_II}, the authors discuss exotic fibers in the context of completed fillings, and a finite number of exotic fibers are present in the completed filling. However, since we are looking at compact fillings, cutting off a cylindrical end, we obtain a disk worth of exotic fibers for every single exotic fiber in the completed filling.} as explained in Remark \ref{rmk: compactified_model}. However, for notational convenience, when we refer to an exotic fiber, we always think of the ``genuinely exotic'' one at $0\in \mathbb C$.

A nearly Lefschetz fibration $\Pi$ is \textbf{allowable} if all vanishing cycles are homologically essential curves on a regular fiber. In analogy with positive allowable Lefschetz fibrations that are studied in the literature and denoted PALFs, we will also refer to these structures as positive allowable nearly Lefschetz fibrations and abbreviate them to PANLFs.

\begin{remark}\label{rem: exotic_nbd_branch}
  Under the biholomorphism $(z_1, z_2) \mapsto (z_2^2 - z_1, z_2)$ from $\C^2 \setminus \{(z^2,z)\}$ to $\C \setminus \{0\} \times \C$, one can identify the neighborhood of an exotic fiber in the above, with the complement of a branch point of a multisection in a Lefschetz fibration as described in \cite[\S 3]{baykur2016multisections}.
\end{remark}

We now sketch the proof of Theorem \ref{thm: Stein} endowing a PANLF a Stein structure. The full proof can be found in \cite{breen-roy-wang_panlf}.

For a bordered Lefschetz fibration, a Weinstein structure can be described by writing down Weinstein models near the regular fibers and the singular fibers, and then patching them up. We write down the details for a manifold $W$ supported by a nearly Lefschetz fibration $\Pi: W\to \Sigma$, with a generic fiber $S$. Both $S$ and $\Sigma$ are compact oriented surfaces with boundary of genera $g_S$ and $g_\Sigma$.
Let $\Exot(\Pi) \subset \Sigma$ denote the images of the (isolated) exotic fibers and let $\Crit(\Pi) \subset \Sigma$ denote the set of (isolated) critical values, which are assumed to be disjoint from $\Exot(\Pi)$.

\begin{figure}
    \centering
    \includegraphics[width=0.45\linewidth]{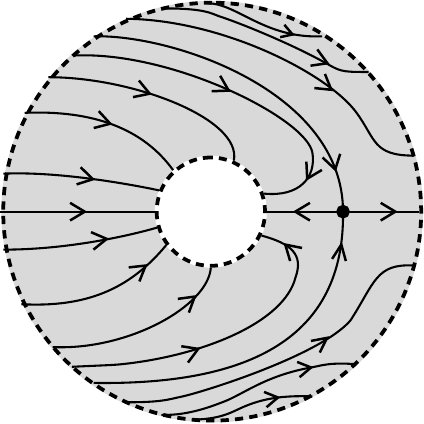}
    \caption{Open annulus and the gradient vector field of $f_A$.}
    \label{fig:annulus}
\end{figure}

As a preparation step, let's first examine a Weinstein model near exotic fibers, since those near the regular and singular fibers have been done previously, e.g. see \cite{ozbagci_survey}. Given $c \in \mathcal{M}_{\exot} \subset \Sigma$, and a sufficiently small neighbourhood $c \in U_c \subset \Sigma$, the preimage $\Pi^{-1}(U_c)$ contains an open subset $V_c \cong A_{\delta} \times \bD^2$ with complex coordinates $(z_1,z_2)$ so that $\Pi$ looks like $z_2^2+z_1$ on $V_c$\footnote{In Definition \ref{def:model}, the local model is given by $z_2^2 - z_1$. We pick the biholomorphic local model $z_2^2+z_1$ here simply so that the vector field $\mathrm{Re}(z_1)$ points to the right.}, with 0 corresponding to $c$. Here $A_{\delta} \coloneqq \{\delta < |z| < 1\}$ denotes an open annulus around 0 in $\C$ as shown in Figure~\ref{fig:annulus}. Let $f_A$ be a function on the annulus whose gradient vector field is shown as in Figure~\ref{fig:annulus}. In particular, $f_A = \text{Re}(z_1)$ for $z_1 \in A_\delta$ with norm close to 1. On $A_{\delta} \times \bD^2$, consider the function $f_{\exot} \coloneqq f_A+x_2^2-y_2^2$.

We can ensure that $f_{\exot}$, along with the standard complex structure on $A_\delta \times \bD^2$ induced from $\C \times \C$, induces a Weinstein structure on $A_\delta \times \bD^2$, compatible with the standard symplectic structure, with the Liouville vector field given by $\nabla f_{\exot}$.

Now, we are ready to describe a particular Weinstein structure on $W$ as follows (see Figure~\ref{fig:panlf_stein}):

\begin{figure}[htbp]
    \vspace{0.4cm}
    \begin{overpic}[width=0.5\textwidth, tics = 10]
    {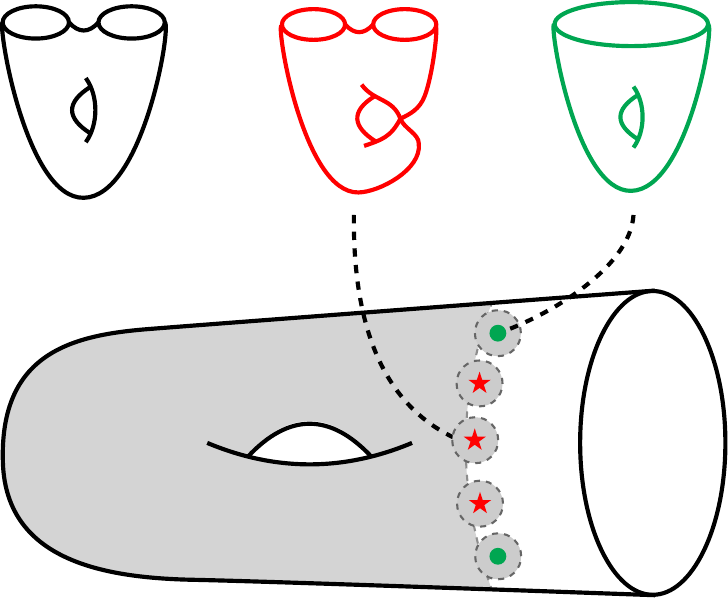}
    \put(3,182){{\scriptsize regular fiber}}
    \put(83,182){{\scriptsize singular fiber}}
    \put(166,182){{\scriptsize exotic fiber}}
    \put(-12,45){{$\Sigma$}}
    \put(65,65){{$\mathcal{U}_1$}}
    \put(160,75){{\scriptsize $U_{p_1}$}}
    \put(153,45){{\scriptsize $U_{q_i}$}}
    \put(157,8){{\scriptsize $U_{p_2}$}}
    \end{overpic}
    \caption{Constructing a Stein structure on a PANLF using local models around critical and exotic points.} %\ar{The open sets around the $U_i$'s should intersect $\mathcal{U}_1$.}  \hm{done}}
    \label{fig:panlf_stein}
\end{figure}

\begin{enumerate}
    \item Consider a Weinstein structure on $S$ given by a function $f_S$ with a single index $0$ critical point and $1-\chi(S)$ index $1$ critical points, and a Weinstein structure on $\Sigma$ given by a function $f_\Sigma$ with a single index $0$ critical point and $1-\chi(\Sigma)$ index $1$ critical points.
    \item Pick a finite open cover with open disks $\mathcal{U} = \cup_x U_x$ of $\Sigma$ such that for every critical value $q \in \Crit(\Pi)$, there is a unique open set $U_q$ containing $q$. Also, we assume that for any image $p \in \Exot(\Pi)$ of an exotic fiber, there is a unique open set $U_p$ containing $p$. We ensure that $\cup_{x \notin \Crit(\Pi)\cup \Exot(\Pi)} U_x =: \cU_1 \cup \cU_2$ is connected, and $\mathcal{U}_1$ is homeomorphic to $\Sigma$.
    \item \textbf{Regular fibers.} Note that $\Pi^{-1}(\mathcal{U}_1)$ is diffeomorphic to $\Sigma \times S$. Consider the function $f_0 = f_S + f_\Sigma$ to obtain a Weinstein structure on this. Since $\cU_1$ is away from the critical points and exotic points, we can ensure that $f_0$ is arbitrarily close to the real part of $\Pi$ in some choice of coordinates over the rest $\mathcal{U}_1 \subset \Sigma$.
    and positively transverse to the boundary of $W$. We can ensure these choices of coordinates agree with the choices made for the local models near critical and exotic fibers.

    \item \textbf{Singular fibers.} Recall that near any critical point $c$ of $\Pi$, there are complex local coordinates $(z_1,z_2)$ on a neighbourhood $V_c \subset W$ such that $\Pi|_{V_c}: V_c \to \Sigma$ looks like $z_1^2+z_2^2$. Consider the function $f_q = x_1^2+x_2^2-y_1^2-y_2^2$ in these coordinates, where $z_j = x_j+iy_j$.
    Given any critical value $q$ of $\Pi$, a critical point $c$ such that $\Pi(c) = q$, and an open set $U_q$ containing $q$, describe $\Pi^{-1}(U_q) = V_c \cup V'_c$ where $V'_c$ is also open and avoids $c$. We can choose $U_q$ so that $U_q \subset \Pi(V_c)$, which in particular ensures that $V'_c$ is a trivial $(S - A)$-bundle over the disk $U_q$, where $A$ is an annulus neighbourhood of the vanishing cycle corresponding to the critical value $q$.
    Since $f_q = \text{Re}(\Pi)$ and $q$ is in the interior, we can argue it matches up with $f_0$ directly on $V'_c \cap \Pi^{-1}(\mathcal{U}_1)$.
    We thus extend $f_0$ to a Weinstein structure on $\Pi^{-1}(\mathcal{U}_1)$ to $\Pi^{-1}(\mathcal{U}_1 \cup \bigcup_{q \in \Crit(\Pi)} U_q)$.
    \item \textbf{Exotic fibers.} Now, near any $p \in \Exot(\Pi)$, we can again find an open set $U_p$ and observe that $\Pi^{-1}(U_p) = V_d \cup V_d'$, where $d \subset V_d$ is a neighbourhood meeting the boundary where $\Pi$ takes the form $\Pi|_{V_d}: \C - B \to \C$, $\Pi(z_1,z_2) = z_2^2+z_1$, with $V_d \cong (\C - B) \times \C$, such that $B$ is a ball around the origin.
    Similarly as above, $V_d'$ is a trivial $S - P$ bundle over the disk, where $P$ is a pair of pants. On $V_d$, define the function $f_p = f_A+x_2^2-y_2^2$ and observe that on $V_d \cap \Pi^{-1}(\mathcal{U}_1)$, the function $f_p$ agrees with $\text{Re}(\Pi)$. Thus we can extend $f_0$ over $\Pi^{-1}(U_p)$, and thus the Weinstein structure is defined on $\Pi^{-1}(\mathcal{U}_1\cup \bigcup_{z \in \Crit(\Pi)\cup\Exot(\Pi)}U_z)$.
    \item Once the above steps are done for every critical and exotic point, the Weinstein structure can be extended with no critical points over the trivial cobordism described by $\Pi^{-1}(\mathcal{U}_2)$.

    \item It will be shown in \cite{breen-roy-wang_panlf} that the Weinstein structure above induces a contact form on the boundary which is a Giroux form for the spinal open book induced on the boundary from the PANLF.
\end{enumerate}

\begin{figure}
    \centering
    \includegraphics[width=0.35\linewidth]{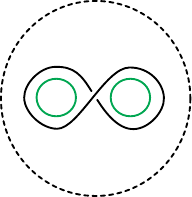}
    \caption{The attaching curve for the Weinstein 2-handle cobordism corresponding to an exotic fiber, living in the neighbourhood of a page of the spinal open book.}
    \label{fig: exotic-cobordism}
\end{figure}

\begin{exercise} \label{ex:panlf_weinstein_handle}
    Show that the boundary interchange mapping class corresponds to a Weinstein cobordism between the relevant contact manifolds, by attaching a Weinstein 2-handle along a Legendrian realization of the knot shown in Figure~\ref{fig: exotic-cobordism}, that lives in a neighbourhood of the pages of the spinal open book. \textit{Hint: Use the local model near the exotic fibers and the definition of the Weinstein function to identify the attaching sphere and the belt-sphere of the handle (up to an isotopy) in neighborhoods of the fibers of $\Pi$, and conclude that the surgery is smoothly a 0-surgery with respect to the framing given by the pages. We give a detailed answer in \cite{breen-roy-wang_panlf}.}
\end{exercise}

The above formulation allows us to write down an explicit Morse function supporting a compatible Weinstein structure on a PANLF. In \cite{breen-roy-wang_panlf} it will be shown that any compatible Weinstein structure on the 4-manifold, i.e., for which the fibers are almost complex, is homotopic to this one.

%%%%%%%%%%%%%%%%%%%%%%%%%%%%%%%%%%%%%%%%%%%%%%%%
\section{Examples of Filling Classifications}\label{sec: examples}
%%%%%%%%%%%%%%%%%%%%%%%%%%%%%%%%%%%%%%%%%%%%%%%%

In this section, we review several known examples and reinterpret them using the language of spinal open books and nearly Lefschetz fibrations to classify their symplectic fillings. Our primary classification tool is Theorem~\ref{thm: monodromy_fact}. In what follows, we outline the strategy for classifying fillings of $(M, \xi)$ supported by a uniform planar spinal open book $(\pi_\Sigma, \pi_P)$, where the page, monodromy, and vertebra are denoted by $(P, \phi, \Sigma = \Sigma_1 \sqcup \dots \sqcup \Sigma_r)$. Note that since we deal with uniform spinal open books, all pages should be diffeomorphic, but the monodromies could be different. Let $m$ be the number of connected components of the paper $M_P$ and $\phi = (\phi_1,\ldots,\phi_m)$. Here, as explained in \S\ref{subsec: factorization}, we assume that we already identify all pages with a single ``model surface'' $P$. The strategy goes as follows:
\begin{enumerate}[(i)]
    \item First, list all possible bases $B \in \mathcal{B}$ such that $(\pi_\Sigma, \pi_P)$ is uniform with respect to $B$. To find such $B$, we need to use Definition~\ref{def:spinalobd} that the number of boundary components of $B$ is the same as the number of connected components of paper $M_P$, namely $m$, and each vertebra is (possibly branched) covering of $B$. These conditions restrict the possible topology of $B$ due to the Riemann--Hurwitz formula:
    \[
        \chi(\Sigma_i) = k_i \cdot \chi(B) - \sum_{p\in \Sigma}(e_p-1),
    \]
    where $k_i$ is the degree of the covering map, which can be directly computed from Definiton~\ref{def:symmetricsobd}, and $e_p$ is the ramification index at $p$, which is equal to $2$ at the branch points and $1$ otherwise, since $\pi_i$ is a simple branched cover.

    \item For each $B \in \cB$, find all positive admissible factorizations of $\phi$ with respect to a monodromy representation $\rho\colon \pi_1(B) \to \Mod(P)$ as follows:

    Suppose each branched covering map $\pi_i: \Sigma_i \to B$ has $n_i$ branch points. Let $n_B = n_1 + \cdots n_r$ be the total number of branch points. This will be the number of exotic fibers in the PANLF induced from the spinal open book and a monodromy factorization.

    Since the spinal open book is uniform, for each boundary component $\partial_j$ of $B$, there exists a corresponding paper component with monodromy $\phi_j$ that intersects the boundary component $\pi_i^{-1}(\partial_j)$ of $\Sigma_i$ for every $i=1,\ldots,r$.

    Now for each $j = 1,\ldots,m$, factorize $\phi_j\circ \rho(\partial_j)^{-1}$ into a product containing boundary interchanges, and positive Dehn twists about homologically essential curves on the page. In fact, each branched cover $\pi_i\colon \Sigma_i \to B$ determines the boundary interchange maps in a factorization. The total number of boundary interchanges across $\phi_1\circ \rho(\partial_1)^{-1},\ldots, \phi_m\circ \rho(\partial_m)^{-1}$ should be $n_B$. Repeat this process for all possible monodromy representations $\rho$.

    \item As a consequence, we obtain a list of symplectic fillings of $(M,\xi)$, each described as a PANLF over a base $B$ with $n_B$ exotic fibers, as $B$ ranges over $\mathcal{B}$. Distinct PANLFs may potentially support equivalent symplectic structures, requiring us to account for double-counting. See Theorem~\ref{thm: parabolicbundles} for an example where different monodromy representations induce equivalent symplectic fillings. However, we still lack an example in which two different bases induce equivalent symplectic fillings.
\end{enumerate}

Note that the symplectic filling corresponding to a factorization above is a PANLF over $B$ with fiber $P$ having $n_B$ exotic fibers and sigular fibers along the curves associated to Dehn twists. Each connected component of paper $M_P$, a $P$-bundle over $S^1$ with monodromy $\phi_j$, corresponds to a subbundle of this PANLF restricted to the boundary component $\partial_j$ of $B$.

In general, the classification problem using the machinery of spinal open books consists of two main parts. First, we seek a planar spinal open book decomposition supporting the contact structure under consideration. Second, we study the monodromy factorization problem in the framed mapping class group. The following examples provide a glimpse of how these procedures are carried out in practice. In many cases, the procedure is much simpler, as the base is a disk (so $\rho$ is trivial), or $P$ has only one component.

%------------------------------------------------
\subsection{The Stein fillable contact structure on the 3-torus} \label{subsec: T3}
%------------------------------------------------
The first example is the Stein fillable contact structure on $\bT^3$, which is arguably the original motivation for spinal open book decompositions. Then, we consider the examples of other torus bundles and finally give some outlook on contact $3$-manifolds supported by more general higher genus open book decompositions.

Recall from \S\ref{subsubsec: T_sobd} that there exists a unique Stein fillable contact structure $\xi_0$ on $\bT^3$. The classification of symplectic fillings for $(\mathbb{T}^3, \xi_0)$ is due to Wendl \cite{Wendl_strongly}.

\begin{theorem}[{\cite{Wendl_strongly}}]
\label{thm: wendl_t3}
    There exists a unique minimal symplectic filling of $(\bT^3, \xi_{0})$ which is given by $\bD^2 \times \bT^2$ equipped with the product symplectic form.
\end{theorem}

Here, we provide an alternate proof of Theorem~\ref{thm: wendl_t3} using spinal open book decompositions, following the strategy outlined above.

\begin{enumerate}[(i)]
    \item In \S\ref{subsubsec: T_sobd}, we construct a spinal open book $\pi = (\pi_{\Sigma},\pi_P)$ supporting $(\bT^3, \xi_0)$ that consists of an annulus page $P$ and two paper components with monodromies $\phi_1=\phi_2=\text{id}$, and two annular vertebrae $\Sigma_1, \Sigma_2$. We first need to list all possible bases $B \in \cB$ such that the spinal open book is uniform with respect to $B$. Recall from Definition~\ref{def:spinalobd} that the number of connected components of $\partial B$ is the same as the number of connected components of the paper $M_P$. Since our $M_P$ has two connected components, we know that $B$ must have two boundary components. According to Definition~\ref{def:spinalobd} again, there exists a (possibly branched) covering map $\pi_i\colon \Sigma_i \to B$ for $i=1,2$. It is straightforward to check $k_1=k_2=1$ from Definition~\ref{def:symmetricsobd}, so $\pi_i\colon \Sigma_i \to B$ are $1$-fold coverings. Thus $B$ must be an annulus. Clearly the covering maps do not have a branch point.  Therefore, we have $\cB=\{A\}$ where $A = S^1 \times I$. In general, $\cB$ can contain more than one element, and we will see such an example in Theorem \ref{thm: example_with_multiple_uniform_surfaces}.

    \item Now, we need to find all possible positive admissible factorizations of $\phi_1$ and $\phi_2$, with respect to $(B,\rho)$, where $\rho\colon \pi_1(B) \to \Mod(P)$ is a monodromy representation. Since $P$ is an annulus, we have $\Mod_\partial(P) = \langle T_{c} \rangle$, where $c$ is the core of $P$. Let $\partial_1$ and $\partial_2$ be the boundary components of $B$ corresponding to the paper components with $\phi_1$ and $\phi_2$, respectively. Since $\partial_2$ is homotopic to $-\partial_1$, we have
    \[
        \rho(\partial_2) = \rho(-\partial_1) = \rho(\partial_1)^{-1}
    \]
    Now, let $\rho(\partial_1) = T_c^k$ for some $k \in \bZ$. Then $\rho(\partial_2) = T_c^{-k}$.

    \begin{exercise}
        Show that $\rho(\partial_1)$ cannot contain a boundary interchange map. \\{\em (Hint: Proposition~\ref{prop:twotwists} may be helpful.)}
    \end{exercise}

    Now we factorize monodromies $\phi_1 \circ \rho(\partial_1)^{-1} = T_c^{-k}$ and $\phi_2 \circ \rho(\partial_2)^{-1} = T_c^{k}$ into a product of positive Dehn twists along homologically essential curves on $P$. There is no boundary interchange since $\pi_1$ and $\pi_2$ do not have a branch point. However, this is impossible unless $k=0$ since otherwise, either $k$ or $-k$ would be negative. Therefore, $k=0$ and the monodromy factorizations of $\phi_1 \circ \rho(\partial_1)^{-1}$ and $\phi_2 \circ \rho(\partial_2)^{-1}$ are both $id$.

    \item $\phi = (\phi_1,\phi_2)$ admits a unique positive admissible monodromy factorization, and hence $(\bT^3, \xi_0)$ admits a unique minimal symplectic filling that admits a bordered Lefschetz fibration over an annular base with an annular fiber without singular and exotic fibers. Therefore, it is diffeomorphic to $A \times A$ and we have
    \begin{align*}
        A \times A &= (S^1 \times I) \times (S^1 \times I)\\
                   &\cong (I \times I) \times (S^1 \times S^1)\\
                   &\cong \bD^2 \times \bT^2
    \end{align*}

\end{enumerate}
    \begin{exercise}
        Show that the product symplectic structure of $\bD^2 \times \bT^2$ is deformation equivalent to the product symplectic structure of $A \times A$.
    \end{exercise}

%---------------------------------------------------
\subsection{Symplectic fillings of torus bundles}
%---------------------------------------------------
The next examples are parabolic torus bundles over $S^1$. Recall from \S\ref{subsec: bundles} that a torus bundle $T_A$ is \textbf{parabolic} if $|\mathrm{tr}(A)|=2$, and it is \textbf{positive} if $A$ is conjugate to $A_k = \begin{pmatrix}
    1 & k\\ 0 & 1
\end{pmatrix}$ for some $k \in \bZ$ and \textbf{negative} if $A$ is conjugate to $-A_k$; there are denoted by $T_{\pm}(k)$. Our main interest here is in rotational contact structures on these bundles, which were studied in \cite{DG:fillings,Van_horn_morris_thesis, ding2018strong, spinal_II,spinal_exotic}. We again follow the strategy outlined at the beginning of this section to classify the symplectic fillings of these contact structures.

\begin{theorem}[Wendl--Lisi--Van Horn-Morris \cite{spinal_II}, Min--Roy--Wang \cite{spinal_exotic}] \label{thm: parabolicbundles}
    Suppose $k \in \bZ$. Then the following statements hold.
    \begin{enumerate}
        \item $(T_+(k),\xi_2)$ is strongly fillable if and only if $k \geq 0$, and for $k \geq 0$, it admits a unique minimal symplectic filling up to deformation equivalence, and this filling is deformation equivalent to a Stein filling.

        \item $(T_-(n), \xi_1)$ is strongly fillable if and only if $k \geq -4$, and for $k \geq -4$, it admits a unique minimal symplectic filling up to deformation equivalence, and this filling is deformation equivalent to a Stein filling.

        \item all other rotational contact structures,  $(T_+(n),\xi_{2k})$ for $k > 1$ and $(T_-(n),\xi_{2j+1})$ for $j \geq 1$, do not admit a strong symplectic filling.
    \end{enumerate}
\end{theorem}

The first part of Theorem~\ref{thm: parabolicbundles} demonstrates that two different monodromy representations, $\rho_1$ and $\rho_2$, can produce equivalent symplectic structures, and the second part provides an example that a factorization of monodromy contains boundary interchange maps. To prove the theorem, we need the following proposition, which provides a factorization of two boundary interchange maps as a product of Dehn twists.

\begin{figure}[htbp]{
  \vspace{0.2cm}
  \begin{overpic}[tics=20]
  {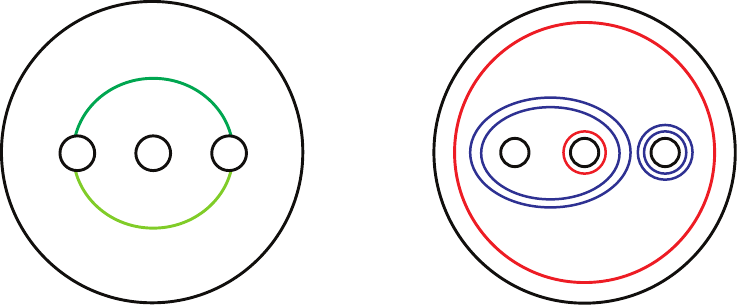}
  \put(70,115){$\alpha$}
  \put(70,24){$\beta$}
  \put(20,134){$1$}
  \put(23,80){$2$}
  \put(58,80){$3$}
  \put(96,80){$4$}
  \end{overpic}}
  \vspace{0.2cm}
  \caption{Left: Two consecutive boundary interchanges $\tau_{\beta}\tau_{\alpha}$ of the inner boundary components along two non-intersecting arcs $\alpha$ and $\beta$. Right: A monodromy factorization of $\tau_{\beta}\tau_{\alpha}$.}
  \label{fig:twotwists}
\end{figure}

\begin{figure}[htbp]{
  \vspace{0.2cm}
  \begin{overpic}[tics=20]
  {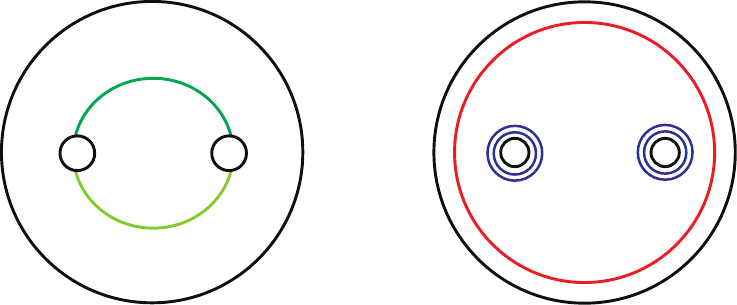}
  \put(70,115){$\alpha$}
  \put(70,24){$\beta$}
  \put(20,134){$1$}
  \put(23,80){$2$}
  \put(96,80){$3$}
  \end{overpic}}
  \vspace{0.2cm}
  \caption{Left: Two consecutive boundary interchanges $\tau_{\beta}\tau_{\alpha}$ of the inner boundary components along two non-intersecting arcs $\alpha$ and $\beta$. Right: A monodromy factorization of $\tau_{\beta}\tau_{\alpha}$.}
  \label{fig:twotwists2}
\end{figure}

\begin{proposition}\label{prop:twotwists}
    Consider a pair of pants with two arcs $\alpha$ and $\beta$ joining inner boundary components as shown in the left drawing of Figure~\ref{fig:twotwists2}. Then make additional boundary component in the middle of the disk as shown in the left drawing of Figure~\ref{fig:twotwists}. Then
    \begin{itemize}
        \item  For a disk with 3 punctures, $\tau_{\beta} \circ \tau_{\alpha} = T_{\partial_1}  T_{\partial_3}  T_{\partial_4}^{-2}  T_{c}^{-2}$, where $c$ is a curve enclosing $\partial_2$ and $\partial_3$. See the right drawing of Figure~\ref{fig:twotwists}.
        \item For a pair of pants, $\tau_{\beta} \circ \tau_{\alpha} = T_{\partial_1} T_{\partial_2}^{-2} T_{\partial_3}^{-2}$. See the right drawing of Figure~\ref{fig:twotwists2}.
    \end{itemize}
\end{proposition}

\begin{exercise}
    Prove Proposition~\ref{prop:twotwists}.
\end{exercise}

\begin{figure}[htbp]{
  \vspace{0.2cm}
  \begin{overpic}[tics=20,width=0.85\textwidth]
  {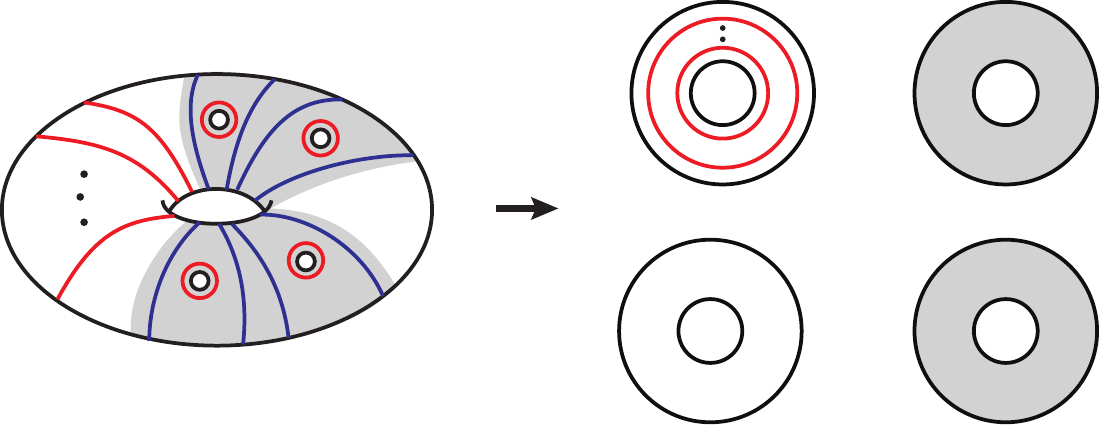}
  \end{overpic}}
  \vspace{0.2cm}
  \caption{Left: An open book decomposition supporting $(T_{+}(k),\xi_{2})$, corresponding to the word $a^k(aba)^{-4}$. Each gray region represents a relative open book with the word $(aba)^{-2}$. Right: A spinal open book decomposition supporting $(T_{+}(k),\xi_{2})$. The gray annuli are vertebrae corresponding to the word $(aba)^{-2}$.}
  \label{fig:sobd_parabolic}
\end{figure}

\begin{proof}[Proof of Theorem~\ref{thm: parabolicbundles}]
    The third part of the theorem is immediate from the fact that the contact structures under consideration contain Giroux torsion. We now consider the first part.

    \begin{enumerate}[(i)]
        \item  We begin by finding a planar spinal open book decomposition supporting $(T_+(k), \xi_2)$. According to Van Horn-Morris \cite[Proof of Theorem~4.3.1]{Van_horn_morris_thesis}, the word for $(T_+(k), \xi_2)$ is $(aba)^{-4}a^k$. Consequently, we obtain a spinal open book with two annular vertebrae, each corresponding to a relative open book with word $(aba)^{-2}$, one annular page with identity monodromy, and another annular page with monodromy $T_c^k$, where $c$ is the core of the annulus (see Figure~\ref{fig:sobd_parabolic}). Thus, we have pages components, both of which are annuli. Here we again denote them by $P$ as a `model page' that identifies two page components. See \S\ref{subsec: factorization} for more details. Also we have monodromy $\phi = (\phi_1, \phi_2)$, where $\phi_1 = T_c^k$ and $\phi_2 = \operatorname{id}$, and vertebrae $\Sigma = \Sigma_1 \sqcup \Sigma_2$, both of which are annuli.

        We then list all possible bases $B \in \cB$ such that the spinal open book is uniform with respect to $B$. Recall from Definition~\ref{def:spinalobd} that the number of connected components of $\partial B$ is the same as the number of connected components of the paper $M_P$. Since our $M_P$ has two connected components, we know that $B$ must have two boundary components. According to Definition~\ref{def:spinalobd} again, there exists a (possibly branched) covering map $\pi_i\colon \Sigma_i \to B$ for $i=1,2$.  It is straightforward to check $k_1=k_2=1$ from Definition~\ref{def:symmetricsobd}, so $\pi_i\colon \Sigma_i \to B$ are $1$-fold coverings. Thus $B$ must be an annulus. Clearly the covering maps do not have a branch point. Therefore, we have $\cB=\{A\}$ where $A = S^1 \times I$.

        \item  Now, we need to find all possible positive admissible factorizations of $\phi_1$ and $\phi_2$, with respect to $(B,\rho)$, where $\rho\colon \pi_1(B) \to \Mod(P)$ is a monodromy representation.  Since $P$ is an annulus, we have $\Mod_\partial(P) = \langle T_{c} \rangle$. Let $\partial_1$ and $\partial_2$ be the boundary components of $B$ corresponding to the paper components with $\phi_1$ and $\phi_2$, respectively. Since $\partial_2$ is homotopic to $-\partial_1$, we have
        \[
            \rho(\partial_2) = \rho(-\partial_1) = \rho(\partial_1)^{-1}
        \]
        Now, let $\rho(\partial_1) = T_c^i$ for some $i \in \bZ$. Then $\rho(\partial_2) = T_c^{-i}$. We factorize monodromies $\phi_1 \circ \rho(\partial_1)^{-1} = T_c^{k-i}$ and $\phi_2 \circ \rho(\partial_2)^{-1} = T_c^{i}$ into a product of positive Dehn twists along homologically essential curves on $P$. There is no boundary interchange since $\pi_1$ and $\pi_2$ do not have a branch point. Now if $k < 0$, this is impossible since either $k-i$ or $i$ would be negative. Therefore, $k\geq 0$ and $k=0$ is already dealt in Theorem~\ref{thm: wendl_t3}. Finally, if $k > 0$, then both $k-i$ and $i$ are positive as long as $0 \leq i \leq k$. Therefore, there are $k+1$ possible monomromy representations $\rho_0, \ldots, \rho_k$ such that $\rho_i(\partial_1) = T^i_c$.

        Then, for each $i = 0,\ldots, k$, we have unique positive admissible factorization
        \[
        \phi_1 \circ \rho_i(\partial_1)^{-1} = T_c^{k-i} \text{\ and \ } \phi_2 \circ \rho_i(\partial_2)^{-1} = T_c^{i}.
        \]

        \item From the factorizations above, we obtain $k+1$ symplectic fillings from PANLF over $B$ without exotic fibers and $k = (k-i) + i$ singular fibers. That is, since $B$ is an annulus, they can be obtained from an honest symplectic fiber bundle over $A$ with fiber $A$ and monodromy representation $\rho_i$ for $i=0,\ldots, k$, and then attaching $k$ Weinstein $2$-handles along the curves associated with the Dehn twists.

        Finally, we need to check for duplications. To do so, we first note that any two oriented $A$-bundles over $A$ are isomorphic to each other, related by a change of trivialization. To see this, consider an $A$-bundle over $A$ with $\rho_0$, which is a trivial bundle $A \times A$. Then the following diffeomorphism induces a bundle isomorphism between two bundles with monodromy representations $\rho_0$ and $\rho_i$:
        \[
            F_i\colon A\times A \to A \times A, \quad F_i(\theta,r,\phi,s) = (\theta,r, \phi - i \theta s,s)
        \]
        That is, after changing trivializations, all bordered Lefschetz fibrations become equivalent and they induce a spinal open book on their boundary with two annular pages, one with monodromy $T_c^k$ and the other with the identity monodromy. Thus $(T_+(k),\xi_2)$ admits a unique minimal symplectic filling if and only if $k\geq0$.
    \end{enumerate}

\begin{figure}[htbp]{
  \vspace{0.2cm}
  \begin{overpic}[tics=20,width=0.85\textwidth]
  {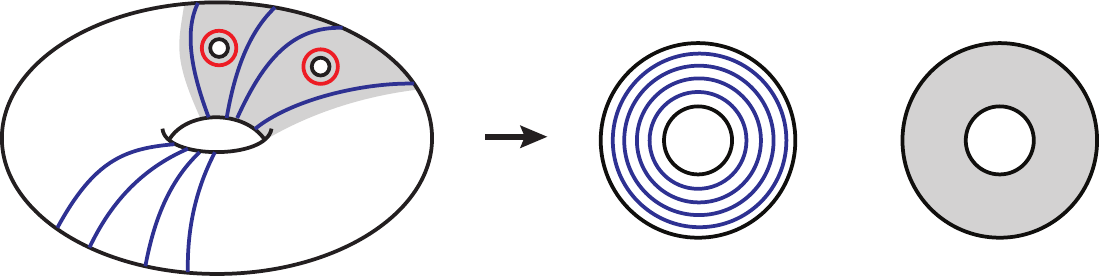}
  \end{overpic}}
  \vspace{0.2cm}
  \caption{Left: An open book decomposition supporting $(T_{-}(-4), \xi_{1})$, corresponding to the word $a^{-4}(aba)^{-2}$. The gray region represents a relative open book with the word $(aba)^{-2}$. Right: A spinal open book decomposition supporting $(T_{-}(-4), \xi_{1})$. The gray annulus is a vertebra corresponding to the word $(aba)^{-2}$.}
  \label{fig:sobd_parabolic_neg}
\end{figure}

    Next, we consider the second part. Again, we begin by finding a spinal open book decomposition supporting $(T_-(k),\xi_1)$.

    \begin{enumerate}[(i)]
        \item According to Van Horn-Morris \cite[Theorem~4.3.1]{Van_horn_morris_thesis}, a word for $(T_-(k), \xi_1)$ is $(aba)^{-2}a^k$. Therefore, we can convert it into a spinal open book decomposition with one annular vertebra $\Sigma$ and one annular page $P$ with monodromy $\phi = T^k_c$ where $c$ is the core of the annular page. See Figure~\ref{fig:sobd_parabolic_neg}.

        We then list all possible bases $B \in \cB$ such that the spinal open book is uniform with respect to $B$. Since our paper $M_P$ has a single connected component, we know that $B$ must have exactly one boundary component. Hence $\chi(B) = 1-2g$ where $g$ is the genus of $B$. It is straightforward to check from Definition~\ref{def:symmetricsobd} that $k=2$. By the Riemann--Hurwitz formula, we obtain
        \[
            0 = 2(1-2g) - b
        \]
        where $b$ denotes the number of branch points. This implies that $g=0$ and $b=2$. Therefore, we have $\cB=\{\bD^2\}$ and $\pi_{\Sigma}\colon \Sigma \to B$ is a $2$-fold branched covering with two branch points.

        \item Now, we need to find all possible positive admissible factorizations of $\phi$, with respect to $(B,\rho)$, where $\rho\colon \pi_1(B) \to \Mod(P)$ is a monodromy representation. Since $B$ is a disk, $\rho$ must be trivial so we only need to factorize $\phi = T_c^k$. Since there are two branch points, a factorization should contain two boundary interchange maps. According to Proposition~\ref{prop:twotwists}, two consecutive boundary interchanges can be factorized into Dehn twists about the boundary components, i.e., $\tau_{\beta} \tau_{\alpha} = T_{\partial_1} T_{\partial_2}^{-2} T_{\partial_3}^{-2}$ as shown in Figure~\ref{fig:twotwists2}. Once we cap off the outer boundary component $\partial_1$, we obtain two consecutive boundary interchanges of an annulus, and its factorization becomes $\tau_{\beta} \tau_{\alpha} = T_{\partial_c}^{-4}$ as $\partial_2$ and $\partial_3$ become homotopic to the core $c$ of the annulus. In our setting, however, two arcs $\alpha$ and $\beta$ may not be homotopic rel boundary, but we still have a Hurwitz equivalent factorization.

        \begin{exercise}
            Let $\alpha$ and $\beta$ be two arcs joining two boundary components of an annulus that are not homotopic rel boundary. Show that $\tau_{\alpha}$ and $\tau_{\beta}$ are related by a framing conjugation (see Definition~\ref{def: generalized_hurwitz_equivalence}).
        \end{exercise}

        Hence we have the following factorization of $\phi$:
        \[
            \phi = T_c^k = T_c^{k+4} \circ \tau_{\beta}\circ \tau_{\alpha}
        \]

        \item  Therefore, $\phi$ admits a unique positive admissible factorization if and only if $k \geq -4$, and this proves the second part of the Theorem.
     \end{enumerate}
\end{proof}

Combined with the results regarding additional tight contact structures on parabolic torus bundles in \cite{spinal_II, Christian:fillings}, Theorem~\ref{thm: parabolicbundles} completes the classification of symplectic fillings of every contact structure on all parabolic torus bundles.

\begin{figure}[htbp]{
  \vspace{0.2cm}
  \begin{overpic}[tics=20]
  {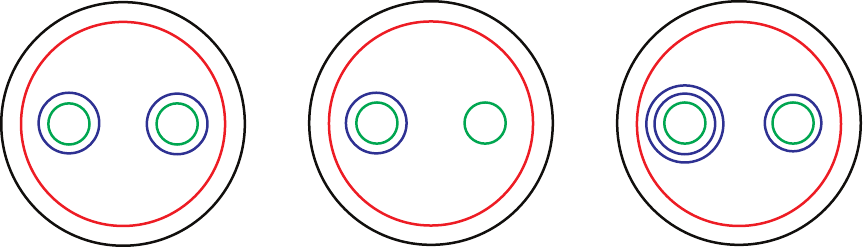}
  \put(55,-12){$\phi_1$}
  \put(205,-12){$\phi_2$}
  \put(355,-12){$\phi_3$}
  \end{overpic}}
  \vspace{0.2cm}
  \caption{The monodromies of the spinal open books supporting the $\pi$-twisting rotational elliptic torus bundles. The green boundary components meet the annulus vertebra.}
  \label{fig:ellipticbundles}
\end{figure}

Next, we consider the elliptic torus bundles with $\pi$-twisting rotational contact structures.

\begin{figure}[htbp]{
  \vspace{0.2cm}
  \begin{overpic}[tics=20]
  {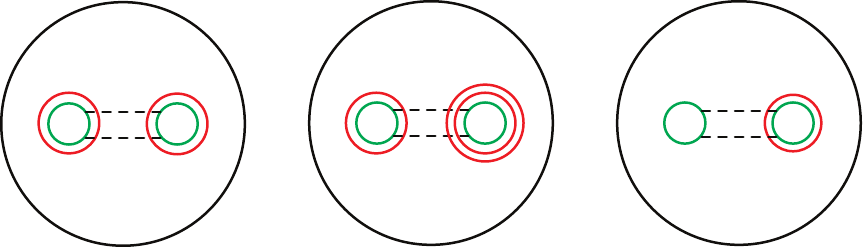}
  \put(55,-12){$\phi_1$}
  \put(205,-12){$\phi_2$}
  \put(355,-12){$\phi_3$}
  \end{overpic}}
  \vspace{0.2cm}
  \caption{Positive admissible monodromy factorizations for $\pi$-twisting elliptic torus bundles, with boundary interchanges along the dotted arcs.}
  \label{fig:ellipticbundles2}
\end{figure}

\begin{theorem}[Min-Roy-Wang \cite{spinal_exotic}]\label{thm: ellipticbundles}
    A $\pi$-twisting rotational contact structure on an elliptic torus bundle admits a unique minimal symplectic filling up to deformation equivalence, and this filling is deformation equivalent to a Stein filling.
\end{theorem}

\begin{exercise}
   Prove Theorem~\ref{thm: ellipticbundles}. \\{\em (Hint: Use Exercise~\ref{exer:elliptic}. See also Figure~\ref{fig:ellipticbundles} and Figure~\ref{fig:ellipticbundles2} for corresponding spinal open books.)}
 \end{exercise}

%-----------------------------------------
\subsection{Further examples}
%-----------------------------------------
Here, we study Stein fillings of some contact 3-manifold supported by a planar spinal open book decomposition that can be uniform with respect to different bases.

This contact 3-manifold is obtained by taking the boundary of the quotient of a product symplectic manifold. We first recall that $\Sigma_{g,n}$ is a genus $g$ surface with $n$ boundary components. Let $\psi_1$ and $\psi_2$ be involutions of $\Sigma_{2,2}$ arising as the deck transformations of the double branched covers $\phi_1: \Sigma_{2,2} \to \Sigma_{1,2}$ and $\phi_2: \Sigma_{2,2} \to \Sigma_{0,2}$, where $\phi_1$ is unbranched and $\phi_2$ has four simple branch points. Also let $\sigma$ be an involution $\sigma(s,t) = (-s,-t)$ of $\Sigma_{0,2} = [-1,1] \times S^1$, with $S^1$ parametrized as $\R / \Z$. Then for $i=1,2$, the following Weinstein domains
\begin{align}\label{eq:E}
  E_i = (\Sigma_{2,2}\times \Sigma_{0,2}) \big/ (z,w) \sim (\psi_i(z), \sigma(w))
\end{align}
induce the same boundary contact 3-manifold $(M,\xi)$, and it is supported by a spinal open book with two annular pages such that each page has a half-twist of the two boundaries as the monodromy, and a single vertebra $\Sigma_{2,2}$. Note that this half-twist is an idemptotent (since $\sigma^2 = id$) and is distinct from the boundary interchange map. We denote the page monodromies by $\phi = (\phi_1,\phi_2) = (\sigma, \sigma)$.%, where $\alpha$ is an arc joining two boundary components of the first page, and $\beta$ is an arc joining two boundary components of the second page.

Although it is challenging to classify all symplectic fillings of $(M,\xi)$ since listing all possible monodromy representations is difficult, we can show that $(M,\xi)$ admits two types of fillings: the first is a bordered Lefschetz fibration (including $E_1$), and the second are nearly Lefschetz fibrations. It is not a priori clear that the symplectic fillings arising as PANLFs in the second class are all symplectic deformation equivalent.

\begin{theorem}[Wendl-Lisi-Van Horn-Morris \cite{spinal_II}, Min-Roy-Wang \cite{spinal_exotic}] \label{thm: example_with_multiple_uniform_surfaces}
    Consider the contact $3$-manifold $(M, \xi)$ above. The Stein fillings of this manifold are either bordered Lefschetz fibrations over $\Sigma_{1,2}$ with $\Sigma_{0,2}$-fiber without singular fibers, or a PANLF over $\Sigma_{0,2}$ with $\Sigma_{0,2}$-fiber, and eight singular fibers.
\end{theorem}

\begin{proof}
    We again follow the strategy outlined at the beginning of the section. Since pages have one boundary component with muliplicity $2$, there exists a (possibly branched) $2$-fold cover $\pi_1\colon \Sigma_{2,2} \to B$ that induces a $2$-fold honest cover on the boundaries. According to the Riemann--Hurwitz formula, $B$ should be either $\Sigma_{1,2}$ or $\Sigma_{0,2}$.

    In the case where $B = \Sigma_{1,2}$, there are no branch points due to the Riemann--Hurwitz formula. Thus there are no exotic fibers, and hence all symplectic fillings are bordered Lefschetz fibrations over $\Sigma_{1,2}$. These can be constructed from a $\Sigma_{0,2}$-bundle over $\Sigma_{1,2}$ with monodromy representation $\rho\colon \Sigma_{1,2} \to \Sigma_{0,2}$, followed by attaching Weinstein $2$-handles along the curves associated with singular fibers that arise from a positive admissible factorization of the page monodromies. We obtain a family of such fillings by varying the monodromy representation.

   In the case where $B = \Sigma_{0,2}$, there are four branch points due to the Riemann--Hurwitz formula. Thus all symplectic fillings contain 4 exotic fibers in their completions, and hence are PANLF over an annular base. These can be constructed from a $\Sigma_{0,2}$-bundle over $\Sigma_{0,2}$ with monodromy representation $\rho\colon \pi_1(\Sigma_{0,2}) \to \Mod(\Sigma_{0,2})$, followed by adding four exotic fibers that contribute boundary interchanges, and attaching Weinstein $2$-handles along curves associated with singular fibers that come from a positive admissible factorization of the page monodromies. There  are four boundary interchanges since there are four branch points for the branched cover map $\Sigma_{2,2} \to \Sigma_{0,2}$.

    The map $\rho$ will be of the form \begin{equation}\label{eq:representation}\rho(c) = \sigma^{\epsilon}T_{c}^k,
    \end{equation} where $\sigma$ represents the involution half-twist in \eqref{eq:E}, $\epsilon \in \{0,1\}$, and $c$ is a based loop around the core of the annulus $\Sigma_{0,2}$.
 
    It follows (similar to the proof of Theorem~\ref{thm: wendl_t3}) that
    \[
    \rho(\partial_1) = T_{c}^k\sigma^{\epsilon}, \rho(\partial_2) = T_c^{-k}\sigma^{-\epsilon}.
    \]

    Since $\phi_1 =\phi_2 = \sigma$, classification of fillings amounts to classifying positive admissible factorizations of
    $\phi_1\circ \rho(\partial_1)^{-1} = T_c^{-k}\sigma^{1+\epsilon}$ and $\phi_2\circ \rho(\partial_2)^{-1} = T_c^k\sigma^{1+\epsilon}$, involving a total of 4 boundary interchanges.

    Let $\tau_b$ represent a boundary interchange on the annulus. We have the following relationship of compactly supported mapping classes (note that on the annulus, despite possibly different choices of arcs, there is a well-defined boundary interchange mapping class and it commutes with $T_c$):

    \[
     \tau_b = T_c^{-2}\sigma, \sigma^2 = id.
    \]

    Now, if $k \neq 0$, we need to find positive admissible factorizations of $T_c^{-k}\sigma^{1+\epsilon}$, and in those, negative Dehn twists can only possibly arise from squares of boundary interchanges as mentioned above, in multiples of 4. It follows (since we can have at most 4 boundary interchanges in any one factorization) that the possible values of $k$ are between 0 and 8. However, if $\epsilon = 0$, we also need to realize a non-compactly supported mapping class using positive admissible factorizations, which means we should have odd powers of boundary interchanges in each factorization.

    \paragraph{Case 1: $0\leq k\leq 4, \epsilon=0$} 3 boundary interchanges for $\partial_1$, 1 for $\partial_2$. $6-k$ singular fibers for $\partial_1$, $k+2$ for $\partial_2$.

    \[
    \phi_1\circ \rho(\partial_1)^{-1} = T_c^{-k}\sigma = T_c^{6-k} \tau_b^3
    \]
    \[
    \phi_2\circ \rho(\partial_2)^{-1} = T_c^k\sigma = T_c^{k+2} \tau_b
    \]
    \paragraph{Case 2: $0\leq k \leq 4, \epsilon=1$} 2 boundary interchanges each for $\partial_1$ and $\partial_2$. $4-k$ singular fibers for $\partial_1$, $k+4$ for $\partial_2$.

     \[
    \phi_1\circ \rho(\partial_1)^{-1} = T_c^{-k}\sigma^2 = T_c^{4-k} \tau_b^2
    \]
    \[
    \phi_2\circ \rho(\partial_2)^{-1} = T_c^k\sigma^2 = T_c^{k+4} \tau_b^2
    \]

    \paragraph{Case 3: $4<k\leq8, \epsilon=0$} Not possible.

    \paragraph{Case 4: $4<k\leq 8, \epsilon=1$} 4 boundary interchanges for $\partial_1$. $8-k$ singular fibers for $\partial_1$, $k$ for $\partial_2$.

    \[
    \phi_1\circ \rho(\partial_1)^{-1} = T_c^{-k}\sigma^2 = T_c^{8-k} \tau_b^4
    \]
    \[
    \phi_2\circ \rho(\partial_2)^{-1} = T_c^k\sigma^2 = T_c^k
    \]

    One can observe in all the cases above, one builds deformation equivalent PANLFs over the annulus. 
\end{proof}

    \begin{exercise}
        Show that the symplectic fillings in the above cases are all deformation equivalent, and are thus deformation equivalent to $E_2$.
    \end{exercise}

\section{Future questions}\label{sec: questions}
In this section, we outline possible future directions of research, and some examples that will be interesting to work out in detail.

%------------------------------------
\subsection{Classifying fillings of higher genus open books}
%------------------------------------
First, we provide a motivating family of examples of open book with higher genus pages. These examples illustrate both the usefulness and complexity of using spinal open books to classify symplectic fillings of contact $3$-manifolds supported by higher genus open books.

Specifically, the following family of open books shows how certain open books with arbitrarily high genus can be converted into planar uniform spinal open books, which makes understanding their fillings a more tractable problem. However these examples also highlight the computational difficulty of the classification problem in this setting.

\begin{figure}[htbp]{
  \vspace{0.2cm}
  \begin{overpic}[tics=20]
  {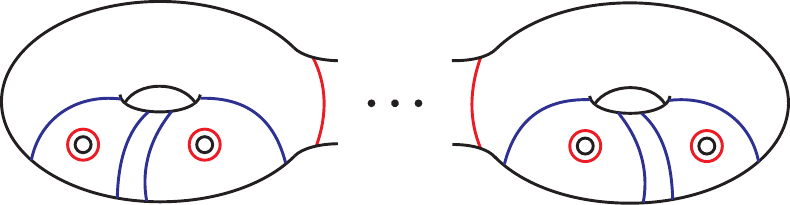}
  \end{overpic}}
  \vspace{0.2cm}
  \caption{A high genus open book decomposition which can be turned into a planar spinal open book. }
  \label{fig:highgenus}
\end{figure}

\begin{example} \label{example: higher genus}
  Consider the family of open books shown in Figure~\ref{fig:highgenus}. For every genus $g\geq 2$, this represents a Stein fillable contact $3$-manifold that admits at least one Stein filling which is a Lefschetz fibration over the disk with genus $g$ fibers, and is at the same time a nearly Lefschetz fibration over the disk with planar fibers. The Stein filling is originally given by \cite{ding2018strong}. In fact, the proof of Proposition 3.2 in \cite{ding2018strong} shows that this open book monodromy admits a positive factorization, which corresponds to a Stein filling that admits a Lefschetz fibration with genus $g$ fibers. On the other hand, we can replace the `handle' parts -- relative open books corresponding to $(aba)^{-2}$ -- with annular vertebrae by Proposition~\ref{prop:obd_ann_vertebra}, and obtain a planar uniform spinal open book supporting the same contact manifold.
  Theorem~\ref{thm: monodromy_fact} then shows that all Stein fillings must have the structure of a PANLF over a disk with planar fibers.
\end{example}

\begin{exercise}
    Use Theorem~\ref{thm: monodromy_fact} to show that the contact $3$-manifold in Example \ref{example: higher genus}, for $g=2$, has at least two minimal symplectic fillings: one with six singular fibers, and one with seven singular fibers.
\end{exercise}

As the solution to the above exercise suggests (the reader can refer \cite{spinal_exotic} for the answer), the filling classification problem deals with complicated configurations of arcs along which boundary interchanges can happen. A resolution to the following question will thus provide handy computational tools in using the technology of planar spinal open books.

\begin{question}
    What is the full classification of minimal symplectic fillings of the contact $3$-manifolds in Example \ref{example: higher genus}?
\end{question}

%------------------------------------------------
\subsection{Giroux correspondence for spinal open books}
%------------------------------------------------
It is shown in \cite[Theorem 1.7]{spinal_I} that the space of Giroux forms supporting a spinal open book is contractible. However, since vertebrae can have complicated topology, it is likely that the notion of \textbf{stabilization} of open books needs to be generalized, for example, including stabilizations on vertebrae.

\begin{question} \label{question: giroux_correspondence}
    What is the correct definition of stabilization of a spinal open book for which the Giroux correspondence holds?
\end{question}

In addition, we need to generalize the notion of right- and left-veering \cite{hkm_rightveering} to the spinal setting.

\begin{question}
    Is there a notion of right- and left-veering, which characterizes overtwistedness for spinal open book decompositions?
\end{question}

We can also ask whether additional symplectic invariants arise that help distinguish fillings beyond topological obstructions, such as the number of vanishing cycles.

\begin{question}
    Are there diffeomorphic but not symplectic deformation equivalent fillings for a fixed contact $3$-manifold?
\end{question}

%------------------------------------------------
\subsection{Support genus problem}
%------------------------------------------------
The first obstruction to the planarity of an open book is given by Theorem \ref{thm: etnyre_negativedefinite}, using McDuff's characterization of ruled surfaces \cite{mcduff_rational_ruled}.

By Theorems \ref{thm: wendl_strongly}, \ref{thm: main} and \ref{thm: Stein}, if a contact $3$-manifold admits a strong but not Stein filling, then it is not supported by a planar (uniform spinal) open book. This is because for planar (uniform spinal) open books, all of its minimal strong fillings are deformation equivalent to positive allowable (nearly) Lefschetz fibrations, which are supported by Stein structures.

Generalizing this, Wendl gives an obstruction to planar uniform spinal open books in ECH.

\begin{theorem}{\cite[Theorem 6]{Wendl_hierarchy}} \label{thm:wendl_U_tower}
    If $(M,\xi)$ is supported by a planar uniform spinal open book, then the ECH contact invariant is at the bottom of a $U$-tower.
\end{theorem}

In light of this, we naturally ask whether a similar obstruction exists in the Heegaard Floer setting.

\begin{question}
    Is there a Heegaard Floer obstruction to a contact three-manifold being supported by a planar spinal open book?
\end{question}

One way to compute the Heegaard Floer contact invariant is by converting a supporting open book $(\Sigma, \phi)$ into a Heegaard diagram. This leads to the following question.

\begin{question} \label{question: Heegaard_diagrams}
    Is there a natural way to convert to spinal open books to Heegaard diagrams?
\end{question}

Recall that there is a description for the relative open book for a spine component with annulus vertebra (see Proposition \ref{prop:obd_ann_vertebra}). One approach to answering Question \ref{question: Heegaard_diagrams} is to find relative open books for spine components with vertebrae of different topology.

 Note that for spinal open books arising as the boundary of PANLFs, there is an algorithm to convert the spinal open book to a regular open book. This is in theory as a result of Theorem \ref{thm: Stein}. In more words, since Exercise~\ref{ex:panlf_weinstein_handle} shows how to produce a Kirby-Weinstein diagram for the manifold, we can convert that into an open book following the algorithm in \cite{avdek_surgery_obd}. However this may not be the most efficient. It is certainly worth asking if there is a direct way to produce regular open books of ``least complexity'' out of spinal open books.

 \begin{question} \label{question: spinal_to_obd}
    What is an efficient way to produce regular open books out of spinal open books?
 \end{question}

One can hope these tools will prove beneficial in answering the following long-standing open question about the support genus of open books. While we know that there are obstructions to contact manifolds being supported by planar open books, we know very little beyond.

\begin{question}[Support genus problem]
    Is there a contact $3$-manifold with support genus two, i.e., there is no genus one open book supporting it?
\end{question}

%--------------------------------------------------
\subsection{Connections to 4-manifold topology and near symplectic structures}
%--------------------------------------------------
Lefschetz fibrations and pencils provide a powerful characterization of symplectic 4-manifolds. Relaxing the geometric condition, Etnyre and Fuller \cite{etnyre-fuller_achiral} showed that an ``achiral Lefschetz fibration'' on a 4-manifold can be endowed with a near-symplectic structure, i.e., a symplectic structure on the complement of a collection of embedded $S^1$'s. They also showed that any smooth, closed, oriented 4-manifold admits an achiral Lefschetz fibration after surgery on a framed circle. Simply put, achiral Lefschetz fibrations allow topological vanishing cycles where the associated monodromy is a left-handed Dehn twist.

One can wonder whether similar notions exist for nearly Lefschetz fibrations, and if they can provide tools for understanding near-symplectic structures. A natural thought is to introduce the clockwise or negative boundary interchange mapping class, to allow for ``topological exotic fibers''.

\begin{question}
    What is the correct notion of an achiral nearly Lefschetz fibration? Does every smooth, closed, oriented 4-manifold admit one? Are they compatible with near symplectic structures?
\end{question}

%%%%%%%%%%%%%%%%%%%%%%%%%%%%%%%

\bibliographystyle{amsalpha}
\bibliography{references}
\end{document}